\documentclass[notitlepage,leqno,10pt]{article}
\usepackage{latexsym}
\usepackage{amsmath}
\usepackage{amsfonts}
\usepackage{amssymb}
\usepackage{amscd}

\renewcommand{\d}{\delta }
\newcommand{\D }{\Delta }

\renewcommand{\l }{\lambda }

\newcommand{\n }{\nabla }

\newcommand{\s }{\sigma }

\renewcommand{\o }{\omega }

\newcommand{\ov}{\overline}
\newcommand{\intbar}{\mathop{\int\makebox(-13.5,0){\rule[4pt]{.7em}{0.3pt}}%
\kern-6pt}\nolimits}

\newcommand{\be}{\begin{equation}}
\newcommand{\ee}{\end{equation}}
\newcommand{\bes}{\begin{equation*}}
\newcommand{\ees}{\end{equation*}}
\newcommand{\ba}{\begin{eqnarray}}
\newcommand{\ea}{\end{eqnarray}}
\newcommand{\bas}{\begin{eqnarray*}}
\newcommand{\eas}{\end{eqnarray*}}
\newenvironment{pf}{\noindent{\sc Proof}.\enspace}{\rule{2mm}{2mm}\medskip}
\newenvironment{pfn}{\noindent{\sc Proof}}{\rule{2mm}{2mm}\medskip}

\newcommand{\R}{\mathbb{R}}

\newcommand{\Z}{\mathbb{Z}}

\newcommand{\N}{\mathbb{N}}
\renewcommand{\o }{\omega }

\author{Martin MAYER$^{a}$,\;\; Cheikh Birahim NDIAYE$^b$}

\date{}

\title{\bf Proof of the remaining cases of the boundary Yamabe problem}

\begin{document}

\newtheorem{lem}{Lemma}[section]
\newtheorem{pro}[lem]{Proposition}
\newtheorem{thm}[lem]{Theorem}
\newtheorem{rem}[lem]{Remark}
\newtheorem{cor}[lem]{Corollary}
\newtheorem{df}[lem]{Definition}

\maketitle

\begin{center}

{\small

\noindent  $^{a, b}$ Mathematisches Institut der Justus-Liebig-Universit\"at Giessen, \\Arndtstrasse 2, D-35392 Giessen, Germany.

}
\
\
{\small

\noindent

}

\end{center}

\footnotetext[1]{E-mail addresses: martin.g.mayer@math.uni-giessen.de, ndiaye@everest.mathematik.uni-tuebingen.de, cheikh.ndiaye@math.uni-giessen.de.}

\

\

\begin{center}
{\bf Abstract}

\end{center}
In this paper, we solve the remaining cases of the boundary Yamabe problem introduced by Escobar\cite{es1} in 1992. Indeed, using the bubbles of Brendle-Chen\cite{bs}, which are an adaptation to manifolds with boundary of the original ones introduced by Brendle\cite{bre2} for the study of the  Yamabe flow on closed Riemannian manifolds of dimension greater or equal to $6$, and the algebraic topological argument of Bahri-Coron\cite{bc}, we solve the cases left open after the work of Brendle-Chen\cite{bs}. Thus combining our work with the ones of Brendle-Chen\cite{bs} and Escobar\cite{es1}, we have that the boundary Yamabe problem is a done deal.

\begin{center}

\bigskip

\centerline{\bf AMS subject classification:  53C21, 35C60, 58J60, 55N10.}

\end{center}
\section{Introduction and statement of the results}
In 1992 Escobar\cite{es1} raised the question of whether every compact \;$n$-dimensional Riemannian manifold with boundary and \;$n\geq 3$\; carries a conformal metric with constant scalar curvature in its interior and vanishing mean curvature on its boundary. In the same work he provides a positive answer unless \;$n\geq 6$, the boundary is umbilic, the manifold is not locally conformally flat and its Weyl tensor vanishes identically on the boundary. Recently, building on the work of Brendle\cite{bre2}, Brendle-Chen\cite{bs} have made a very important progress on the cases left open by Escobar\cite{es1} by solving many situations and reducing the remaining ones to the positivity of the ADM mass of a certain class of asymptotically flat manifolds.  Unfortunately, the latter positivity {\em is not know to hold}.
\vspace{6pt}

\noindent
Our main goal in this work is to solve the cases left open by Brendle-Chen\cite{bs} and Escobar\cite{es1}. Indeed using the bubbles of Brendle-Chen\cite{bs} combined with a suitable scheme of the algebraic topological argument of Bahri-Coron\cite{bc}, we give a positive answer to the cases of the boundary Yamabe problem remaining after the works Brendle-Chen\cite{bs} and Escobar\cite{es1}, by showing a result which covers all the cases left open after the latter works. To state clearly our theorem, we first fix some notation. Given \;$(\ov M, g)$\; a compact \;$n$-dimensional Riemannian manifold with boundary $\partial M$, interior $M$, and \;$n\geq 3$,  we denote by \;$L_g=-4\frac{n-1}{n-2}\D_g+R_g$\; the conformal Laplacian of $(\ov M, g)$ and \;$B_g=\frac{4(n-1)}{n-2}\frac{\partial}{\partial n_g}+2(n-1)H_g$\; the conformal Neumann operator of $(M, g)$, with \;$R_g$\; denoting the scalar curvature of \;$(\ov M, g)$, \;$\D_g$ \;denoting the Laplace-Beltrami operator with respect to $g$, \;$H_g$\; is the mean curvature of \;$\partial M$\; in \;$(\ov M, g)$,  \;$\frac{\partial}{\partial n_g}$ is the outer Neumann operator on $\partial M$ with respect to $g$. Furthermore, we define the following boundary Yamabe functional 
\begin{equation}\label{eq:Yamabefunctional}
\mathcal{E}_g(u):=\frac{\langle L_gu, u\rangle+\langle B_g u, u\rangle}{(\int_{M}u^{\frac{2n}{n-2}}dV_g)^{\frac{n-2}{n}}}, \;\;\;\;\;u\in W^{1, 2}_{+}(\ov M),
\end{equation}
where \;$\langle L_gu, u\rangle:=\langle L_gu, u\rangle_{L^2(M)}$, $\langle B_gu, u\rangle:=\langle B_gu, u\rangle_{L^2(\partial M)}$, $dV_g$ is the volume form with respect to $g$, and $W^{1, 2}_+(\ov M):=\{u\in W^{1, 2}(\ov M): \;\;u>0\}$\; with $W^{1, 2}(\ov M)$\; denoting the usual Sobolev space of functions which are $L^2$-integrable with their first derivatives (for more information, see \cite{aubin} and \cite{gt}). Moreover, we recall that the Yamabe invariant of $(M, \partial M, g)$ is defined as 
\begin{equation}\label{eq:Yamabeinvariant}
 \mathcal{Y}(M, \partial M, g):=\inf_{u\in W^{1, 2}_+(\ov M)}\mathcal{E}_g(u).
\end{equation}
Now, having fixed the needed notation, we are ready to state our theorem which reads as follows.
\begin{thm}\label{eq:existence}
Assuming that \;$(\ov M, g)$\; is a $n$-dimensional compact Riemannian manifold with boundary \;$\partial M$\; and interior \;$M$\; such that \;$\partial M$\; is umbilic in \;$(\ov M, g)$, \;$n\geq 6$, and \;$\mathcal{Y}(M, \partial M, g)>0$,\; then \;$\ov M$\; admits a Riemannian metric conformal to \;$g$\; with constant scalar curvature in \;$M$\; and zero mean curvature on \;$\partial M$. 
\end{thm}
Hence as already said, Theorem \ref{eq:existence} and the works of Brendle-Chen\cite{bs} and Escobar\cite{es1} imply the following positive answer to the Yamabe problem for manifolds with boundary introduced by Escobar\cite{es1}.
\begin{thm}\label{eq:BoundaryYamabe}
Every \;$n$-dimensional compact Riemannian manifold with boundary and \;$n\geq 3$\; carries a conformal metric with constant scalar curvature in the interior and zero mean curvature on the boundary.
\end{thm}
\begin{rem}\label{eq:Maybemin}
We would like to emphasize the fact that we do not know whether our solution is a minimizer of the Euler-Lagrange functional \;$\mathcal{E}_g$\; or not.
\end{rem}
\vspace{10pt}

\noindent
Under the assumptions of Theorem \ref{eq:existence}, the boundary Yamabe problem is equivalent to finding a smooth and positive solution to the following semilinear elliptic boundary value problem\\
\begin{equation}\label{eq:bvp}
\left\{
\begin{split}
L_gu&=u^{\frac{n+2}{n-2}} \;\;&\text{in}\;\;M,\\
 B_gu&=0\;\;&\text{on}\;\;\partial M.
\
\end{split}
\right.
\end{equation} 
On the other hand, thanks to the work of Cherrier\cite{cherrier}, we have that smooth solutions of \eqref{eq:bvp} can be found by looking at critical points of the boundary Yamabe functional \;$\mathcal{E}_g$, and we will pursue such an approach here. Precisely, to prove Theorem \ref{eq:existence}, we use the algebraic topological argument of Bahri-Coron\cite{bc}. The latter argument is well known to be related to the notion of critical points at infinity introduced by Bahri\cite{bah}. However, to study the critical points at infinity of a conformally invariant geometric variational problem (precisely which are so on model cases after "blow-ups"), the first issue is to find optimal variational bubbles which are variational bubbles, namely bubbles which reduce the variational analysis of the geometric problem to the one of a functional defined on the space of parameters responsible of the conformal invariance, and which are optimal with respect to the approach one wants to use to handle the variational analysis of the latter reduced functional. The algebraic topological argument or sometimes called the barycenter technique of Bahri-Coron\cite{bc} can be implemented  with Morse-Conley theory as done by the the second author for the \;$Q$-curvature problem. But, such an approach requires a complete understanding of the critical points at infinity of the variational problem, and hence requires optimal variational bubbles where optimality is with respect to the establishment of a Morse lemma at infinity. On the other hand, the algebraic topological argument of A. Bahri and J. Coron, as they have shown in \cite{bc}, can be applied without the study of the critical points at infinity and with bubbles which are not necessarily optimal variational bubbles where optimality is with respect to the derivation of a  Morse lemma at infinity, but need variational bubbles which are sharp in a certain sense as observed by Bahri-Brezis\cite{BB}, and we refer to them as algebraic topological bubbles. On the other hand, in the field of existence of solutions of conformally invariant geometric variational problems and the long time convergence of their associated geometric flows, some type of bubbles have appeared, starting from the work of Aubin\cite{au}, to the best of our knowledge (see also the work of Hebey-Vaugon\cite{hv}), and we call those compactifying bubbles, because of their purpose of use. In 1984 Schoen\cite{schoen} proposes some sharp compactifying bubbles to treat the cases left open after the work of Aubin\cite{au}. Recently, in his study of the asymptotic behaviour of the Yamabe flow on closed Riemannian manifolds of dimension greater or equal to \;$6$, Brendle\cite{bre2} introduced some bubbles with a very interesting geometry, which allows him to reduce the convergence of the Yamabe flow on closed Riemannian manifolds to the positivity of the ADM mass of a certain class of asymptotically flat manifolds. Thus one can speculate that Brendle\cite{bre2}'s bubbles are optimal compactifying bubbles, where optimal means here best for showing convergence of the Yamabe flow on closed Riemannian manifolds or showing the existence of a minimizer of the Yamabe functional. Furthermore, we observe from the geometry behind Brendle\cite{bre2}'s bubbles that they solve the same variational problem as optimal variational bubbles to first order, where optimality is with respect to the establishment of a Morse Lemma at infinity. Hence, with the conclusion on the Yamabe flow obtained by Brendle\cite{bre2}, one can even speculate that they are optimal variational bubbles in the sense of Morse theory, an hence are algebraic topological bubbles. In this work, we make the latter heuristic rigorous for the boundary adaptation of the Brendle\cite{bre2}'s bubbles, namely the ones of Brendle-Chen\cite{bs}. Indeed, we first use the bubbles of Brendle-Chen\cite{bs} and some estimates of them established in \cite{bs}, to show that they are sharp enough to bring us to an ``entrance door'' at infinity sufficiently close to the critical points at infinity 
of \;$\mathcal{E}_g$\;, where closeness is measured with respect to the discrepancy in energy from the critical values at infinity and verification of the interaction identity at infinity. Next, we use the latter information to show that they are algebraic topological bubbles by running a suitable scheme of the barycenter technique of A. Bahri and J. Coron for the case under study and deduce Theorem \ref{eq:existence}.
\vspace{6pt}

\noindent
Some comments are in order regarding our strategy of proof. Even if our argument for proving Theorem \ref{eq:existence} is based on an application of the algebraic topological argument of Bahri-Coron\cite{bc}, it differs at the level of presentation and at the technical level with previous works using such a method. Indeed, our presentation follows more closely the variant of the algebraic topological argument of Bahri-Coron\cite{bc} developed by the second author\cite{nd6} for the supercritical $Q$-curvature problem on closed Riemannian manifolds. Furthermore and more importantly, it differs from the other applications of the barycenter technique of Bahri-Coron\cite{bc} in the fact that those works use optimal variational bubbles in the sense of Morse theory which is possible because of being in low dimensions or in special geometry, while here we rely on Brendle-Chen\cite{bs}'s bubbles for which we do not know whether they are optimal variational bubbles in the sense of establishing a Morse lemma at infinity and showing that turns out to {\em not be a trivial matter}. Hence, at the technical level, we are armed in a weaker way, since the optimal variational bubbles used in previous applications of the algebraic topological argument of Bahri-Coron\cite{bc} verify the sharp estimate that they differ from the {\em standard} bubbles in uniform topology to order $O\left(\frac{1}{\l^{n-2}}\right)$ where $\l$ is the concentration parameter (see \cite{bah1}, \cite{BB}, \cite{bc}, \cite{gam1}, \cite{gam2}, \cite{ould1}), while for Brendle-Chen\cite{bs}`s bubbles, it is known that such a property {\em does not hold} as easily seen from the work of Khuri-Marques-Schoen\cite{kms} (see also the work of Disconzi-Khuri\cite{diskhuri}).
\vspace{6pt}

\noindent
\begin{rem}\label{eq:maybemimimizer}
We would like to point out that the method of this paper is used in our forthcoming paper \cite{martndia} in combination with the Chen\cite{chen}'s bubbles to settle the remaining cases of the Riemann mapping problem of Escobar\cite{es2}. Furthermore, it is also used in our forthcoming paper \cite{martndia1} in combination with the techniques of Brendle\cite{bre2}, Brendle-Chen\cite{bs}, and  Chen\cite{chen} to settle the remaining cases of the problem of existence of conformal metrics with constant scalar curvature and constant mean curvature studied by Araujo\cite{araujo1},\cite{araujo2}, Escobar\cite{es3} and Han-Li\cite{hanli1},\cite{hanli2}.
\end{rem}
\vspace{6pt}
\begin{rem}\label{eq:masscpi}
We would like to make a remark about the work of Brendle-Chen\cite{bs} and ours. In order to better make our point, we first recall that to the best of the authors knowledge, there are only three approaches to solve the Yamabe problem on Riemannian manifolds with or without boundary. They are the minimization technique (see \cite{au}, \cite{bs},  \cite{chen}, \cite{es1}, \cite{es2}, \cite{schoen}), the flow approach (see \cite{bre1}, \cite{bre2}, \cite{struwe}), and the algebraic topological method (\cite{BB}, \cite{ould1}). All these approaches rely on the quantization phenomenon and are linked via the critical points at infinity of the associated variational problem. Furthermore, the first one is based on positivity of mass distribution, while the second one needs on top of that the strong interaction phenomenon. However, like the flow approach, the algebraic topological method needs quantization phenomenon and strong interaction phenomenon, but it requires only boundedness of mass distribution. Thus all these methods are linked to general relativity via ADM mass as shown by Schoen\cite{schoen}, Brendle\cite{bre1}, \cite{bre2}, Ben-Chen-Chtioui-Hammami\cite{bcch}, Brendle-Chen\cite{bs},  Chen\cite{chen}, and it is not surprising that Brendle-Chen\cite{bs} could not conclude, while here we can. In fact their work indicates highly the possibility to use their bubbles to perform a suitable algebraic topological argument and indeed this is the case as we show in this paper.
\end{rem}
\begin{rem}
We would like to point out the fact that, we believe that the Brendle\cite{bre1}'s bubbles (dimension greater or equal to \;$6$) and Schoen\cite{schoen}'s bubbles (dimension less or equal than \;$5$\; and locally conformally flat case) are optimal variational bubbles, in the sense that they can be used to establish a  Morse Lemma at infinity for the Yamabe functional via a singularity analysis of the Yamabe flow and using them as parametrization for the shadow flow, as it is being investigated in \cite{martndia3}.  Such a Morse Lemma at infinity will shed some light on the Weyl vanishing conjecture and the ADM mass of the class of asymptotically flat manifolds introduced by Brendle\cite{bre2} and Schoen\cite{schoen}. In fact, it will imply that the Weyl vanishing conjecture is true along the Yamabe flow and is a fundamental identity at infinity for the  critical points at infinity of the Yamabe functional, meaning that the flow is compact in the region violating it. Furthermore, it will imply that the ADM mass of the class of asymptotically flat manifolds introduced by Brendle\cite{bre2} and Schoen\cite{schoen} are the coefficients of the one-bubble ends in the concentration parameters of the critical points at infinity of the Yamabe functional. On the other hand, a scheme to find optimal variational bubbles for Yamabe-type functionals is proposed in \cite{martndia3} to settle the Bahri-Brezis conjecture on nonlinear elliptic equations with critical Sobolev growth on closed Riemannian manifolds.
\end{rem}
\vspace{8pt}

\noindent
The structure of this paper is as follows. In Section \ref{eq:notpre}, we fix some notation and give some preliminaries, like the set of formal barycenters of \;$\partial M$\; and present some useful topological properties of them. Furthermore, we recall the Brendle-Chen\cite{bs}'s bubbles and the fact that they can be used to replace the standard bubbles in the analysis of diverging Palais-Smale (PS) sequences of the Euler-Lagrange functional \;$\mathcal{E}_g$. Moreover, using a result of Brendle\cite{bre2} and another one of Brendle-Chen\cite{bc}, we derive two estimates for them which permit to enter in the world of algebraic topological bubbles. In Section \ref{eq:mappingbary}, we carry the technical preparation needed to show that Brendle-Chen\cite{bs}'s bubbles are algebraic topological ones, by using our estimates to map the space of barycenter of \;$\partial M$\; of any order into suitable sublevels of \;$\mathcal{E}_g$\; via the Brendle-Chen\cite{bs}'s bubbles in such a way that a suitable scheme of the barycenter techniques of A. Bahri and J. Coron is readily applicable. Finally, in Section \ref{eq:algtop}, we define the neighborhood of potential critical points at infinity of $\mathcal{E}_g$ and present the latter readily applicable algebraic topological argument to prove Theorem \ref{eq:existence}.
\vspace{8pt}

\noindent
\begin{center}
{\bf Acknowledgements}
\end{center}
C. B. Ndiaye has been supported  by the DFG project "Fourth-order uniformization type theorems for $4$-dimensional Riemannian manifolds".
\section{Notation and preliminaries}\label{eq:notpre}
In this section, we fix some notations and give some preliminaries. First of all, from now until the end of the paper \;$(\ov M, g)$\; will be the given underlying compact $n$-dimensional Riemannian manifold with boundary \;$\partial M$\; and interior \;$M$, $\partial M$\; is totally geodesic in \;$(\ov M, g)$, $n\geq 6$, and \;$\mathcal{Y}(M, \partial M, g)>0$. 
\vspace{6pt}

\noindent
In the following, for any Riemannian metric \;$\bar g$\; on \;$\ov M$, we will use the notation \;$B^{\bar g}_{p}(r)$\; to denote the geodesic ball with respect to $\bar g$ of radius \;$r$\;and center \;$p$. We also denote by \;$d_{\bar g}(x,y)$\; the geodesic distance with respect to $\bar g$ between two points \;$x$\;and \;$y$\; of \;$\ov M$. $inj_{\bar g}(\ov M)$\;stands for the injectivity radius of \;$(\ov M, \bar g)$, $dV_{\bar g}$\;denotes the Riemannian measure associated to the metric\;$\bar g$, and $dS_{\bar g}$\; the volume form on \;$\partial M$\; with respect to the metric induced by \;$\bar g$\; on \;$\partial M$. For simplicity, we  use \;$B_p(r)$\; to denote \;$B^g_{p}(r)$, namely\; $B_p(r)=B^g_p(r)$. For \;$a\in \ov M$,  we use the notation \;$\exp_a^{\bar g}$\; to denote the exponential map with respect to \;$\bar g$\; and set for simplicity $\exp_a:=\exp_a^g$.
\vspace{6pt}

\noindent
$\N$\;denotes the set of nonnegative integers, $\N^*$\;stands for the set of positive integers, and  for $n\in \N^*$, $\R^n$ the standard $n$-dimensional Euclidean space, $\R^n_+$ the open positive half-space of $\R^n$, and $\bar \R^n_+$ its closure in $\R^n$. For simplicity, we will use the notation \;$\R_+:=\R^1_+$, and $\bar \R_+:=\bar \R^1_+$. For $r>0$, $B^{\R^n}_0(r)$ denotes the open ball of \;$\R^n$\; of center \;$0$\; and radius \;$r$, $\bar  B^{\R^n}_0(r)$ for its closure in \;$\R^n$, $B^{\R^n_+}_0(r):=B^{\R^n}_0(r)\cap\bar \R^n_+$ , and $\bar B^{\R^n_+}_0(r):=\bar B^{\R^n}_0(r)\cap\bar \R^n_+$. $S^n$ denotes the unit sphere of \;$\R^{n+1}$, and \;$S^n_+$ \;the positive spherical cap, namely $S^n_+:=S^n\cap \R^{n+1}_+$. We use \;$g_{S^n}$\; to denote the round metric on $S^n$ and $g_{S^n_+}$\; its restriction on \;$S^n_+$. For \;$p\in \N^*$,  $\ov M^p$, $M^p$, and \;$(\partial M)^p$ denote respectively the cartesian product of $p$ copies of $\ov M$,  \;$M$ and \;$\partial M$. For \;$p\in \N^*$, $F((\partial M)^p)$ denotes the fat diagonal of $(\partial M)^p$, namely \;$F((\partial M)^p):=\{A:=(a_1, \cdots, a_p)\in ( \partial M)^P:\;\;\exists\;i\neq j\;\;\text{with}\;a_i=a_j\}$.  We define \;$((\partial M)^2)^*:=(\partial M)^2\setminus Diag((\partial M)^2)$ where $Diag((\partial M)^2)$ is the diagonal of \;$(\partial M)^2$, namely $Diag((\partial M)^2):=\{(a, a): \;\,a\in \partial M\}$. For $p\in \N^*$, $\sigma_p$ stands for the permutations group of \;$p$\; elements and \;$\D_{p-1}$ the following simplex \;$\D_{p-1}:=\{(\alpha_1, \cdots, \alpha_p):\alpha_i\geq 0, i=1, \cdots, p, \sum_{i=1}^p\alpha_i=1\}$.
\vspace{6pt}

\noindent
For $1\leq p\leq \infty$ and $k\in \N$, $\beta\in  ]0, 1[$, $L^p(M)$\; and \;$L^p(\partial M)$, \;$W^{k, p}(M)$, $C^k(\ov M)$, and \;$C^{k, \beta}(\ov M)$ stand respectively for the standard $p$-Lebesgue space on $M$ and $\partial M$, $(k, p)$-Sobolev space, $k$-continuously differentiable space and $k$-continuously differential space of H\"older exponent \;$\beta$, all with respect to \;$g$ (if the definition needs a metric structure) and for precise definitions and properties, see  for example \cite{aubin} or \cite{gt}.
\vspace{6pt}

\noindent
 For  \;$a\in \ov M$, $O_a(1)$ stands for quantities bounded uniformly in $a$.  For $\epsilon$ positive and small, and \;$a\in \ov M$, $O_{a, \epsilon}(1)$ stands for quantities uniformly bounded in $a$ and $\epsilon$. For $\epsilon$ positive and small, $o_{\epsilon}(1)$ means quantities which tend to \;$0$\; as \;$\epsilon$\; tends to $0$. For $\l$ large and \;$a\in \ov M$, \;$O_{a, \l}(1)$ stands for quantities uniformly bounded in $a$ and $\l$. For \;$a\in \ov M$, \;$\epsilon$\; and $\delta$ positive and small, and $\l$ large, $O_{a, \epsilon, \delta}(1)$ and \;$O_{a, \l}(1)$\; stand respectively for quantities which are bounded uniformly in \;$a$, $\epsilon$, and\;$\d$, and in $a$ and \;$\l$. For $a\in \ov M$, $\epsilon$ positive and small, and $\l$ large, $o_{a, \epsilon}(1)$ and $o_{a, \l}(1)$ stand respectively for quantities which tend to \;$0$\; uniformly in $a$ as $\epsilon$ tends to $0$, and as  $\l$ tends to $+\infty$. For \;$A\in (\partial M)^2$ and $\l$ large, $O_{A, \l}(1)$ and $o_{A, \l}(1)$ stands respectively for quantities which are bounded uniformly in \;$A$\; and \;$\l$, and which tend to $0$ uniformly in $A$ as $\l$ tends to $+\infty$. For \;$p\in \N^*$, $A\in (\partial M)^p$, $\bar \alpha\in \D_{p-1}$, and \;$\l$ large, $O_{A, \bar \alpha, \l}(1)$\; and \;$o_{A, \bar \alpha, \l}(1)$ stand respectively for quantities which are uniformly bounded in \;$p$, $A$, $\bar\alpha$, and $\l$ and for quantities which tend to $0$ uniformly in \;$p$, $A$, and \;$\bar\alpha$ \;as \;$\l$\; tend to $+\infty$. For $x\in \R$, we will use the notation \;$O(x)$\; and $o(x)$ to mean  respectively\;$|x|O(1)$\; and $|x|o(1)$ where \;$O(1)$\; and \;$o(1)$\; will be specified in all the contexts where they are used. Large positive constants are usually denoted by \;$C$\; and the value of \;$C$\; is allowed to vary from formula to formula and also within the same line. Similarly small positive constants are also denoted by \;$c$\; and their values may vary from formula to formula and also within the same line.   The symbol \;$\sum_{i\neq j}$\;always means a double sum over the associated index set under the assumption \;$i\neq j$.                                                                                                                                                                                                                                                                                                                                                                                                                                                                                                                                                                                                                                                                                                                                                                                                                                                                                                                                                                                                                                                                                                                                     \vspace{6pt}

\noindent
For \;$X$ a topological space, \;$H_{*}(X)$\; will denote the singular homology of \;$X$ with \;$\Z_2$\; coefficients, and $H^*(X)$ for the cohomology. For \;$Y$ a subspace of \;$X$, $H_*(X, Y)$ will stand for the relative homology. The symbol \;$\frown$\; will denote the cap product between cohomology and homology. For a map $f:X\rightarrow Y$, with \;$X$\; and \;$Y$\; topological spaces, \;$f_*$\; stands for the induced map in homology, and \;$f^*$\; for the induced map in cohomology.
For $p\in \N$, we set
\begin{equation}\label{eq:defenergylevel}
W_p:=\{u\in W^{1, 2}_+(\ov M):\;\;\;\;\mathcal{E}_g(u)\leq (p+1)^{\frac{n}{2}}\mathcal{Y}(S^n_+)\},
\end{equation}
where
\begin{equation}\label{eq:yamabehsphere}
\mathcal{Y}(S^n_+):=\mathcal{Y}(S^n_+, \partial S^n_+, g_{S^n_+}).
\end{equation}
\vspace{6pt}

\noindent
For a Riemannian metric $\bar g$ defined on $\ov M$, we denote by \;$G_{\bar g}$\; the Green's function of $(L_{\bar g}, B_ {\bar g})$ satisfying the normalization
\begin{equation}\label{eq:greennorm}
\lim_{d_{\bar g}(a, x) \longrightarrow 0}(d_{\bar g}(a, x))^{n-2}G_{\bar g}(a, x)=1,
\end{equation} 
and set
\begin{equation}\label{eq:greenback}
G:=G_g.
\end{equation}
On the other hand, using the existence of conformal normal coordinates (see \cite{gun}, \cite{lp}, and \cite{marques}) and recalling that  \;$\partial M$\; is totally geodesic in $(\ov M, g)$, we have that for every large positive integer\;$m$\; and for every \; $a \in \ov M$, there exists a positive function $u_a\in C^{\infty}(\ov M)$ such that the metric $g_a = u_a^{\frac{4}{n-2}}g$ verifies
\begin{equation}\label{eq:detga}
det g_a(x)=1 +O_{a, x}((d_{g_a}(a, x))^m)\;\;\text{for}\;\;\; x\in B^{g_a}_a( \varrho_a),
\end{equation}
with \;$O_{a, x}(1)$\; meaning bounded by a constant independent of \;$a$\; and \;$x$, $0<\varrho_a<\frac{inj_{g_a}(\ov M)}{10}$. Moreover, we can take the family \; $u_a$, $g_a$, and $\varrho_a$\; such that
\begin{equation}\label{eq:varro0}
\text{the maps}\;\;\;a\longrightarrow u_a, \;g_a\;\;\text{are}\;\;C^0\;\;\;\text{and}\;\;\;\;\frac{1}{4}\geq\varrho_a\geq \varrho_0>0,
\end{equation}
for some small positive \;$\varrho_0$\; satisfying \;$\varrho_0<\frac{inj_g(\ov M)}{10}$, and
\begin{equation}\label{eq:proua}
\begin{split}
&||u_a||_{C^2(\ov M)}=O_a(1),\;\;\frac{1}{\ov C^2} g\leq g_a\leq \ov C^2 g, \;\;\;a\in \ov M, \\\;&u_a(x)= 1+O_a(d^2_{g_a}(a, x))=1+O_a(d_{g}^2(a, x)) \;\;\text{for}\;\;x\in\;\;B_a^{g_a}(\varrho_0)\supset B_a(\frac{\varrho_0}{2\ov C}),\;\;\;\;a\in \ov M,\\&
u_a(a)=1,\;a\in \ov M\;\;\;R_{g_a}(a)=0,\;\;a\in  M\;\;\text{and}\;\; H_{g_a}=0\;\;a\in \partial M,
\end{split}
\end{equation}
for some large positive constant \;$\ov C$\; independent of \;$a$, and for the meaning of \;$O_a(1)$ in \eqref{eq:proua}, see section \ref{eq:notpre}. For \;$a\in \ov M$, and \;$\epsilon$\; positive, 
we define the standard bubbles as follows
\begin{equation}\label{eq:standardbubbles}
\d_{a, \epsilon}(x):=\left(\frac{\epsilon}{\epsilon^2+(d_{g_a}(a, x))^2}\right)^{\frac{n-2}{2}}.
\end{equation}
For $a\in \ov M$ and $0<r<\varrho_0$, we set also 
\begin{equation}\label{eq:greena}
G_a:=G_{g_a}, \;\;\;\;\exp_a^a=\exp_{a}^{g_a}\;\;\;\;\text{and}\;\;\;B_a^a(r):=B_a^{g_{a}}(r).             
\end{equation}
On the other hand, the conformal invariance properties of the couple conformal Laplacian and conformal Neumann operator implies
\begin{equation}\label{eq:invarpro}
\begin{split}
&\mathcal{E}_g(u)=\mathcal{E}_{g_{a}}(u_a^{-1}u), \;\;\;\;\int_Mu^{\frac{2n}{n-2}}dV_g=\int_{M}(u_a^{-1}u)^{\frac{2n}{n-2}}dV_{g_a}\;\;\;\;\;\text{for}\;\;\;\;u\in W^{1, 2}_+(M),\\ &G_g(x, y)=G_{a}(x, y)u_a(x)u_a(y),\;\; \;\;(x, y)\in \ov M^2\;\;\text{and}\;\;\;\;a\in \ov M.
\end{split}
\end{equation}
We also define the following quantities
\begin{equation}\label{eq:defc}
c_0:=4n(n-1), \;\;\;c_1:=\int_{\R^n_+}\left(\frac{1}{1+|x|^2}\right)^ndx,  \;\;\text{and}\;\;\;\;c_2:=4\frac{n-1}{n-2}\int_{\R^n_+}\left|\n\left[\left(\frac{1}{1+|x|^2}\right)^{\frac{n-2}{2}}\right]\right|^2dx.
\end{equation}
Furthermore, we set
\begin{equation}\label{defc3}
c_3:=\int_{\R^n_+}\left(\frac{1}{1+|x|^2}\right)^{\frac{n+2}{2}}dx,
\end{equation}
and define the following quantity which depends only on $(\ov M, g)$
\begin{equation}\label{eq:defcg}
c_g=\frac{c_3}{4c_1}\min_{((\partial M)^2)^*}G,
\end{equation}
and for the definition of $((\partial M)^2)^*$, $G$ and $c_3$, see above. We recall that the numbers $c_i$ ($i=0, 1, 2$) and $\mathcal{Y}(S^n_+)$ verify the following relation 
\begin{equation}\label{eq:relationcy}
 c_2=c_0c_1\;\;\;\text{and}\;\;\;\mathcal{Y}(S^n_+)=\frac{c_2}{c_1^{\frac{n-2}{n}}}.
\end{equation}
Moreover, for \;$p\in \N^*$, $A:=(a_1, \cdots, a_p)\in (\ov M)^p$, \;$\bar \l:=(\l_1, \cdots, \l_p)\in (\R_+)^p$, we associate the following quantities (which naturally appear in the analysis of diverging PS sequences of the Euler-Lagrange functional \;$\mathcal{E}_g$)
 \begin{equation}\label{eq:varepsilonij}
     \varepsilon_{i, j}:=\varepsilon_{i, j}(A, \bar \l):=\frac{c_{3}}{(\frac{\l_i}{\l_j}+\frac{\l_j}{\l_i}+\lambda_i\l_jG^{\frac{2}{2-n}}(a_{i},a_{i}))^{\frac{n-2}{2}}}, \;\;i, j=1, \cdots, p, \;\;i\neq j.                                                                                                                                                                                                                                                                                                                                                                                                                                                                                                                                                                                                                                                                                                                                                                                                                                                                                                                                                                                                                                                                                                                                                                                           \end{equation}
                                                                  
\vspace{6pt}

\noindent
Now, we are going to present some topological properties of the space of formal barycenter of $\partial M $ that we will need for our algebraic topological argument for existence. To do that, for $p\in \N^*$, we recall that the set of formal barycenters of \;$\partial M$\; of order $p$ is defined as follows, 
 \begin{equation}\label{eq:barytop}
B_{p}(\partial M):=\{\sum_{i=1}^{p}\alpha_i\d_{a_i}: \;\;a_i\in \partial M, \;\alpha_i\geq 0,\;\; i=1,\cdots, p,\;\,\sum_{i=1}^{p}\alpha_i=1\}, 
\end{equation}
and set
\begin{equation}
B_0(\partial M)=\emptyset. 
\end{equation}
Furthermore, we have the existence of \;$\Z_2$\; orientation classes \;$w_p\in H_{np-1}(B_{p}(\partial M), B_{p-1}(\partial M))$\; and that the cap product acts as follows 
\begin{equation}\label{eq:actioncap}
\begin{CD}
 H^l(B_p(\partial M)\setminus B_{p-1}(\partial M))\times H_k(B_{p}(\partial M), B_{p-1}(\partial M))@>\frown>> H_{k-l}(B_{p}(\partial M), B_{p-1}(\partial M)).
 \end{CD}
\end{equation}
Moreover, there holds
\begin{equation}\label{eq:purem}
B_p(\partial M)\setminus B_{p-1}(\partial M)\simeq ((\partial M)^p)^*\times_{\sigma_p}\D_{p-1}^*,
\end{equation}
where
$$
((\partial M)^p)^*:=(\partial M)^p)\setminus F((\partial M)^p),
$$
and
$$
\D_{p-1}^*:=\{(\alpha_1, \cdots, \alpha_p)\in \D_{p-1}:\;\; \;\alpha_i>0, \;\;\forall\;i=1, \cdots, p\}.
$$
On the other hand, since $\partial M$ is a closed $(n-1)$-dimensional manifold, then we have 
\begin{equation}\label{eq:defom}
\text{an orientation class}\,\;0\neq O^{*}_{\partial M}\in H^{n-1}(\partial M).
\end{equation}
Furthermore, there is a natural way to inject $\partial M$ into $((\partial M)^p)^*\times_{\sigma_p}\D_{p-1}^*$, namely an injection
\begin{equation}\label{eq:inj}
i: \partial M\longrightarrow ((\partial M)^p)^*\times_{\sigma_p}\D_{p-1}^*,
\end{equation}
such that
\begin{equation}\label{eq:defom1}
i^*(O^*_p)=O^*_{\partial M}\;\; \text{with}\;\;0\neq O^{*}_p\in H^{n-1}(((\partial M)^p)^*\times_{\sigma_p}\D_{p-1}^*).
\end{equation}
Identifying $O^*_{\partial M}$ and $O^*_p$ via  \eqref{eq:defom1}, and using \eqref{eq:actioncap} and \eqref{eq:purem}, we have the following well-know formula, see \cite{kk}.
\begin{lem}\label{eq:transfert}
There holds 
$$
\begin{CD}
 H^{n-1}(((\partial M)^p)^*\times_{\sigma_p}\D_{p-1}^*)\times H_{np-1}(B_{p}(\partial M), B_{p-1}(\partial M))@>\frown>> H_{np-n}(B_{p}(\partial M), B_{p-1}(\partial M))\\@>\partial>>H_{np-n-1}(B_{p-1}(\partial M), B_{p-2}(\partial M)),
 \end{CD}
$$
and 
\begin{equation}\label{eq:classt}
\omega_{p-1}=\partial (O^*_{\partial M}\frown w_p).
\end{equation}
\end{lem}
\vspace{8pt}

\noindent
Next, we are going to discuss some important properties of the Brendle-Chen\cite{bs}'s bubbles. In order to do that, we first recall that, in his study of the longtime behaviour of the Yamabe flow on closed Riemannian manifolds of dimension greater or equal to \;$6$, Brendle\cite{bre2} has introduced a very interesting family of bubbles to replace the standard bubbles in the analysis of diverging Palais-Smale (PS) sequences of the Yamabe functional and which verify a sharp Yamabe-energy estimate. Later, jointly with S. Chen, in \cite {bs}, they extend the techniques of Brendle\cite{bre2} to the case of Riemannian manifolds with boundary and define a family of bubbles centered at boundary points and with similar properties. Indeed, for $\delta$ small and positive, they define a family of bubbles $v_{a, \epsilon, \delta}$ (see formula (4.2) in \cite{bs}), $a\in \partial M$ and $\epsilon$ positive and small such that they can replace the standard bubbles in the analysis of diverging PS sequences of $\mathcal{E}_g$. Precisely, $v_{a, \epsilon, \delta}$ is defined as a suitable perturbation of the standard bubbles glued with an appropriate scale of the Green's function \;$G_a$\; centered at $a$ as follows
\begin{equation}\label{eq:bsbubbles}
\begin{split}
v_{a, \epsilon, \delta}(\cdot)
=
\chi_{\delta}(\cdot)(\d_{a, \epsilon}(\cdot)+w_{a, \epsilon}(\cdot))
+
(1-\chi_{\delta}(\cdot))\epsilon^{\frac{n-2}{2}}G_{a}(a, \cdot),
\end{split}\end{equation}
and
$$
\frac{\partial v_{a, \epsilon, \delta}}{\partial n_{g_a}}=0\;\;\text{on}\;\;\partial M,
$$
where 
\begin{equation}\label{eq:cutoff}
\chi_{\delta}(x):=\chi\left(\frac{d_{g_a}(a, x)}{\delta}\right),
\end{equation}
and \;$\chi$ is a cut-off function defined on \;$\bar \R_+$\; satisfying \;$\chi$ is non-negative, $\chi(t)=1$ if $t\leq 1$ and $\chi(t)=0$ if $t\geq 2$, $\delta_{a, \epsilon}$ is defined as in \eqref{eq:standardbubbles}, $G_{a}(a, \cdot)$ is defined as in \eqref{eq:greena}, and in normal coordinates around $a$ with respect to $g_a$, we have that \;$w_{a, \epsilon}$\; satisfies the following pointwise estimate
\begin{equation}\label{eq:wdecay}
\vert \partial^{\beta}w_{a, \epsilon}(x)\vert
\leq C_n(|\beta|)
\frac
{\epsilon^{\frac{n-2}{2}}}
{(\epsilon^{2}+r^{2})^{\frac{n-4+\beta}{2}}}\;\;\text{with}\;\;\;r=d_{g_a}(a, x)\; \;\text{and}\; \;\;x\in B_a^a(\varrho_0),
\end{equation}
where $\varrho_0$ is as in \eqref{eq:proua} and $C_n(|\beta|)$ is a large positive constant which depends only on $n$ and $|\beta|$. Furthermore, and more importantly, they verify the following energy estimate which is a weak form of Proposition 4.1 in  \cite{bre1}, but sufficient for the purpose of this paper.
\begin{lem}\label{eq:brenchenenergy}
There exists \;$0<\d_0\leq\varrho_0$ small such that for every \;$0<2\epsilon\leq\delta\leq \d_0$ and for every $a\in \partial M$, there holds
\begin{equation}\label{eq:brenchenenergyest}
\begin{split}
\int_M\left(4\frac{n-1}{n-2}|\n_{g_a} v_{a, \epsilon, \delta}|^2+R_{g_a}v_{a, \epsilon, \delta}^2\right)dV_{g_a}\leq \mathcal{Y}(S^n_+)\left(\int_Mv_{a, \epsilon, \delta}^{\frac{2n}{n-2}}dV_{g_a}\right)^{\frac{n-2}{n}}-\epsilon^{n-2}\mathcal{I}(a, \delta)\\+O_{a, \epsilon, \delta}(\d^2\epsilon^{n-2}+\epsilon^n\delta^{-n}),
\end{split}
\end{equation}
where \;$\mathcal{Y}(S^n_+)$\; is defined by \eqref{eq:yamabehsphere}, \;$\mathcal{I}(a, \delta)$\; is a flux integral verifying \;$\mathcal{I}(a, \delta)=O_{a, \delta}\left(1\right)$, and for the meaning of\; $O_{a, \delta}\left(1\right)$\; and \;$O_{a, \epsilon, \delta}\left(1\right)$, see Section \ref{eq:notpre}.
\end{lem}
\vspace{6pt}

\noindent
On the other hand, recalling that the Brendle-Chen\cite{bs}'s bubbles are boundary adaptation of the Brendle\cite{bre2}'s bubbles and using the same argument as the one of Proposition 20 in \cite{bre2}, we have that they verify the following interaction estimates.
\begin{lem}\label{eq:brencheninteract}
There exists a large constant $C_1>0$ such that for every $2\epsilon_1\leq 2\epsilon_2\leq\delta\leq \delta_0$ and every $a_1, a_2\in \partial M$, there holds
 \begin{equation}\label{eq:brencheninteractest}
 \begin{split}
 \int_M v_{a_1,  \epsilon_1, \delta}\left|-4\frac{n-1}{n-2}\D_{g_{a_2}} v_{a_2, \epsilon_2, \delta}+R_gv_{a_2, \epsilon_2, \delta}-c_0v_{a_2, \epsilon_2,\delta}^{\frac{n+2}{n-2}}\right|dV_{g_{a_2}}\leq\\ C_1\left(\delta^2+\frac{\epsilon_2^2}{\delta^2}\right)\left(\frac{\epsilon_2^2+d_{g_{a_2}}^2(a_1, a_2)}{\epsilon_1\epsilon_2}\right)^{\frac{2-n}{2}},
 \end{split}
\end{equation}
where $c_0$ is defined by \eqref{eq:defc}.
\end{lem}
\vspace{6pt}

\noindent
Furthermore, using \eqref{eq:bsbubbles}, it is easy to see that the following estimate holds.
\begin{lem}\label{eq:bsvolest}
Assuming that \;$0<\epsilon \leq\delta_0^{\frac{n}{2}}$\; and \;$a\in \partial M$, then we have
\begin{equation}\label{eq:estconformalvol}
\int_Mv_{a, \epsilon, \epsilon^{\frac{2}{n}}}^{\frac{2n}{n-2}}dV_{g_a}=c_1+o_{a, \epsilon}(1),
\end{equation}
where $c_1$ is as in \eqref{eq:defc}, and for the meaning of $o_{a, \epsilon}(1)$, see Section \ref{eq:notpre}.
\end{lem}
\vspace{6pt}

\noindent
Thus, setting, 
\begin{equation}\label{eq:val}
 v_{a}^{ \l}:=v_{a, \frac{1}{\l}, (\frac{1}{\l})^{\frac{2}{n}}}, \;\;a\in \partial M,\;\;\;\;\l\geq \frac{2}{\delta_0^{\frac{n}{2}}}
\end{equation}
and
\begin{equation}\label{eq:varphialb}
\varphi_{a, \l}:=u_av_{a}^{\l}, \;\;a\in \partial M, \;\;\;\;\l\geq \frac{2}{\delta_0^{\frac{n}{2}}},
\end{equation}
where \;$\delta_0$ is still given by Lemma \ref{eq:brenchenenergy}, we have clearly that Lemma \ref{eq:brenchenenergy}, Lemma \ref{eq:brencheninteract}, and Lemma \ref{eq:bsvolest} combined with \eqref{eq:invarpro} imply the following lemmata which will play an important role in our application of the barycenter technique of Bahri-Coron\cite{bc}.
\begin{lem}\label{eq:bubbleestl}
Assuming that $a\in \partial M$ and \;$\l\geq \frac{2}{\delta_0^{\frac{n}{2}}}$, then the following estimate holds
\begin{equation}\label{eq:bubblesestenergy}
\mathcal{E}_g(\varphi_{a, \l})\leq \mathcal{Y}(S^n_+)\left(1+O_{a, \l}\left(\frac{1}{\l^{n-2}}\right)\right),
\end{equation}
where $\mathcal{Y}(S^n_+)$ is as in \eqref{eq:yamabehsphere} and for the meaning of \;$O_{a, \l}\left(1\right)$, see Section \ref{eq:notpre}.
\end{lem}
\begin{lem}\label{eq:bubbleinteractl}
 There exists a large constant $C_2>0$ such that for every $a_1, a_2\in \partial M$, and for every $\l\geq \frac{2}{\delta_0^{\frac{n}{2}}}$, we have
 \begin{equation}\label{eq:bubbleinteractest}
 \begin{split}
 \int_M \varphi_{a_1,  \l}\left|-4\frac{n-1}{n-2}\D_g\varphi_{a_2, \l}+R_g\varphi_{a_2, \l}-c_0\varphi_{a_2, \l}^{\frac{n+2}{n-2}}\right|dV_g\leq\\ C_2\left[\left(\frac{1}{\l}\right)^{\frac{4}{n}}+\left(\frac{1}{\l}\right)^{\frac{2n-4}{n}}\right]\left(1+\l^2d_g^2(a_1, a_2)\right)^{\frac{2-n}{2}},
 \end{split}
 \end{equation}
 where $c_0$ is as in \eqref{eq:defc}.
\end{lem}\begin{lem}\label{eq:bsvolestl}
Assuming that  $a\in \partial M$ and $\l\geq \frac{2}{\delta_0^{\frac{n}{2}}}$, then there holds
\begin{equation}\label{eq:estconformalvol}
\int_M \varphi_{a, \l}^{\frac{2n}{n-2}} dV_{g}=c_1+o_{a, \l}(1),
\end{equation}
where $c_1$ is as in \eqref{eq:defc} and for the meaning of \;$o_{a, \l}(1)$, see Section \ref{eq:notpre}.
\end{lem}

\vspace{6pt}
\noindent
On the other hand, using \eqref{eq:bsbubbles}-\eqref{eq:wdecay}, and \eqref{eq:val}, we have that \;$v_{a}^{\l}$\; decomposes as follows 
\begin{equation}\label{eq:decompval}
v_{a}^{ \l}(\cdot)=\chi^{\l}(\cdot)\left(\delta_{a}^{ \l}(\cdot)+w_{a}^{\l}(\cdot)\right)\left(1-\chi_{\l}(\cdot)\right)\frac{G_{a}(a, \cdot)}{\l^\frac{{n-2}}{2}}
\end{equation}
where
\begin{equation}\label{eq:wal}
w_{a}^{ \l}:=w_{a, \frac{1}{\l}},\;\;\;\delta_{a}^{\l}:=\delta_{a, \frac{1}{\l}}, \;\;\text{and}\;\;\;\chi^{\l}=\chi_{(\frac {1}{\l})^{\frac{2}{n}}},
\end{equation}
and \;$w_{a}^{\l}$\; satisfies the following pointwise estimate
\begin{equation}\label{eq:walest}
\vert \partial^{\beta}w_{a}^{\l}(x)\vert
\leq C_n(|\beta|)
\frac
{\l^{\frac{n-6}{2}+\beta}}
{(1+\l^2r^{2})^{\frac{n-4+\beta}{2}}}\;\;\text{with}\;\;\;r=d_{g_a}(a, x)\; \;\text{and}\; \;\;x\in B_a^a(\varrho_0),
\end{equation}
Now, for $p\in \N^*$,  and $A:=(a_1, \cdots, a_p)\in (\partial M)^p$ and $\l\geq \frac{2}{\delta_0^{\frac{n}{2}}}$, we associate the following quantities
\begin{equation}\label{eq:intweight}
\epsilon_{i, j}:=\epsilon_{i, j}(A, \l):= \int_M \varphi_{a_i, \l}^{\frac{n+2}{n-2}}\varphi_{a_j, \l}dV_g, \;\;\;i, j=1, \cdots, p, \;\;i\neq j.
\end{equation}
Using \eqref{eq:varphialb}, \eqref{eq:decompval}-\eqref{eq:intweight}, we have the following lemma which provides self and interaction estimates, and a relation between \;$\epsilon_{i, j}(A, \l)$\; and \;$\varepsilon_{i, j}(A, \bar \l)$ with $\bar \l:=(\l, \cdots, \l)$, and for the meaning of \;$\varepsilon_{i, j}(A, \bar \l)$ see \eqref{eq:varepsilonij}.
\begin{lem}\label{eq:interactestl}
Assuming that $p \in \N^*$, $A:=(a_1,\cdots, a_p)\in (\partial M)^p$ and $\l\geq \frac{2}{\d_0^{\frac{n}{2}}}$, then \\
1) For every $i, j=1, \cdots, p$ with $i\neq j$, we have \\
i) 
\begin{equation*}
\begin{split}
\epsilon_{i,j}\longrightarrow 0
\Longleftrightarrow
\varepsilon_{i,j}\longrightarrow 0,
\end{split}
\end{equation*}
where $\varepsilon_{i, j}:=\varepsilon_{i, j}(A, \bar  \l)$ with $\bar \l:=(\l, \cdots, \l)$ and $\epsilon_{i, j}:=\epsilon_{i, j}(A, \l)$, and for their definitions see respectively \eqref{eq:varepsilonij} and \eqref{eq:intweight}.\\
ii)\\
There exists $\;0<C_3<\infty$\; independent of $p$, $A$ and $\l$ such that the following estimate holds
\begin{equation}
\begin{split}
C_3^{-1}<\frac{\epsilon_{i,j}}{\varepsilon_{i,j}}<C_3.
\end{split}
\end{equation} 
iii)\\
If \;$\varepsilon_{i,j} \longrightarrow 0$, then 
\begin{equation*}
\begin{split}
\epsilon_{i,j}=(1+o_{ \varepsilon_{i, j}}(1))\varepsilon_{i,j},
\end{split}
\end{equation*}
and for the meaning of $o_{\varepsilon_{i, j}}(1)$, see Section \ref{eq:notpre}.\\
2) For every \;$i=1, \cdots, p$, there holds
\begin{equation*}
\begin{split}
\langle L_{g}\varphi_{a_i, \l}, \varphi_{a_i, \l}\rangle+\langle B_g \varphi_{a_i, \l}, \varphi_{a_i, \l}\rangle=c_0(1+o_{a
_i, \l}(1))\int_M \varphi^{\frac{2n}{n-2}}_{a_i, \l}dV_g,
\end{split}
\end{equation*}
where $c_0$ is given by \eqref{eq:defc} and for the meaning of $o_{a_i, \l}(1)$, see Section \ref{eq:notpre}.\\
3) For every \;$i, j=1,\cdots, p$ with $i\neq j$, there holds
\begin{equation*}
\begin{split}
\langle L_{g}\varphi_{a_i, \l}, \varphi_{a_j, \l}\rangle+\langle B_g\varphi_{a_i, \l}, \varphi_{a_j, \l}\rangle=(1+o_{A_{i, j}, \l}(1))c_{0}\epsilon_{i,j},\;\;\;\;\text{and}\;\;\;\;\;
\epsilon_{j,i}=(1+o_{A_{i, j}, \l}(1))\epsilon_{i,j},
\end{split}
\end{equation*}
where $A_{i, j}:=(a_i, a_j)$ and for the meaning of $o_{A_{i, j}, \l}(1)$, see Section \ref{eq:notpre}.
\end{lem}
\begin{pf}
 First of all, to simplify notation, for every  $i=1, \cdots, p$\; we set
 \begin{equation}\label{eq:auxilary}
 \begin{split}
&\varphi_{i}:=\varphi_{a_i, \l}, \;\; v_i:=v_{a_i}^{ \l},\;\;\d_i:=\delta_{a_i}^{\l}, \;\; w_{i}:=w_{a_i}^{ \l}, \;\;G_i(\cdot):=G_{a_i}(a_i, \cdot),\\& B_i:=B_{a_i}^{a_i}(\d_0), \;\;\;\exp_i:=\exp_{a_i}^{a_i}\; \;\text{and}\;\;u_i:=u_{a_i},
\end{split}
\end{equation}
and for the meaning of $B_{a_i}^{a_i}(\delta_0)$ and $\exp_{a_i}^{a_i}$, see Section \ref{eq:notpre}.
Next, using  \eqref{eq:walest} and \eqref{eq:auxilary}, we have that $v_i$ verifies the following pointwise estimate
\begin{equation}\begin{split}
v_{i}
= &
\left(1+o_{a_i, \l}(1)\right)\left(\chi^{\l}\d_i+(1-\chi^{\l})\frac{G_i}{\l^{\frac{n-2}{2}}}\right).
\end{split}\end{equation}
Moreover, using \eqref{eq:greennorm}, \eqref{eq:standardbubbles}, \eqref{eq:wal}, and \eqref{eq:auxilary}, we obtain
\begin{equation}\label{eq:estp1}
\begin{split}
\chi^{\l}\d_{i}
= &
\chi^{\l}\left(\frac{\l}{1+\l^2r^{2}}\right)^{\frac{n-2}{2}}
=
\chi^{\l}\left(\frac{\l}{1+\l^2G_{i}^{\frac{2}{2-n}}\frac{r^{2}}{G_{i}^{\frac{2}{2-n}}}}\right)^{\frac{n-2}{2}} \\
= &
\left(1+o_{a_i, \l}(1)\right)\chi^{\l}\left(\frac{\l}{1+\l^2G_{i}^{\frac{2}{2-n}}}\right)^{\frac{n-2}{2}}
\end{split}\end{equation}
where $r$ is as in \eqref{eq:wdecay} with $a$ replaced by $a_i$, and
\begin{equation}\label{eq:estp2}
\begin{split}
\left(1-\chi^{\l}\right)\left(\frac{\l}{1+\l^2G_{i}^{\frac{2}{2-n}}}\right)^{\frac{n-2}{2}}
= &
\left(1-\chi^{\l}\right)\frac{G_{i}}{\l^{\frac{n-2}{2}}}\left(\frac{1}{1+O_{a_i, \l}(\vert \frac{1}{\l}\vert^{\frac{2(n-2)}{n}})}\right) \\
= &
\left(1+o_{a_i, \l}(1)\right)\left(1-\chi^{\l}\right)\frac{G_i}{\l^{\frac{n-2}{2}}}.
\end{split}\end{equation}
Hence, combining  \eqref{eq:estp1} and \eqref{eq:estp2}, we obtain
\begin{equation}\label{eq:viest}
\begin{split}
v_{i}=\left(1+o_{a_i, \l}(1)\right)\left(\frac{\l}{1+\l^2G_{i}^{\frac{2}{2-n}}}\right)^{\frac{n-2}{2}}.
\end{split}\end{equation}
Now, using \eqref{eq:proua}, \eqref{eq:varphialb}, \eqref{eq:intweight}, \eqref{eq:auxilary}, and \eqref{eq:viest}, we derive the following estimate for $\epsilon_{i, j}$ ($i, j=1, \cdots, p$ and $i\neq j$)
\begin{equation}\label{eq:41}
\begin{split}
\epsilon_{i,j}
= &
\int_M v_{i}^{\frac{n+2}{n-2}}v_{j}\frac{u_{j}}{u_{i}}dV_{g_{a_{i}}} \\
= &
\int_{B_i}v_{i}^{\frac{n+2}{n-2}}v_{j}\frac{u_{j}}{u_{i}}dV_{g_{a_{i}}}
+
O_{A_{i, j}, \l}\left(\frac{1}{\l^n}\right)\\=&\left(1+o_{A_{i, j}, \l}(1)\right)
u_{j}(a_{i})\int_{B^{\l}(0)}\left(\frac{1}{1+|x|^{2}}\right)^{\frac{n+2}{2}}
\left[(\frac{1}{1+\l^{2}G_{j}^{\frac{2}{2-n}}(\exp_{i}(\frac{x}{\l}))})^{\frac{n-2}{2}}\right] dx\\&+
O_{A_{i, j}, \l}\left(\frac{1}{\l^n}\right),
\end{split}\end{equation}
where 
\begin{equation}
B^{\l}(0):=B^{\R^n_+}_0(\l \delta_0),      
\end{equation}
and for the meaning of $B^{\R^n_+}_0(\l \d_0)$, see Section \ref{eq:notpre}. From \eqref{eq:41} it follows that 
\begin{equation}\label{eq:equivalence0}
\begin{split}
\epsilon_{i,j}\longrightarrow 0
\Longleftrightarrow
\l^2 G_{j}^{\frac{2}{2-n}}(a_{i})\longrightarrow+ \infty
\Longleftrightarrow
\varepsilon_{i,j}\longrightarrow 0.
\end{split}\end{equation}
Thus, we have that the proof of i) of point 1) is complete. Now, since $\epsilon_{i,j}$ and $\varepsilon_{i,j}$ are bounded by definition, then thanks to \eqref{eq:equivalence0}, to prove ii) of point 1), we can assume without loss of generality that
\begin{equation}\label{eq:large1}
\begin{split}
\l^{2}G_{j}^{\frac{2}{2-n}}(a_{j})\gg 1. 
\end{split}
\end{equation} 
Thus under the latter assumption, setting
\begin{equation}
\mathcal{A}=B^{\R^n_+}_{0}(\gamma\l\sqrt{G_{j}^{\frac{2}{2-n}}(a_{i})})
\end{equation} 
for $\gamma>0$ small and using Taylor expansion, we obtain that the following estimate holds on\; $\mathcal{A}$
\begin{equation}\label{eq:45}
\begin{split}
&\left(\frac{1}{1+\l^{2}G_{j}^{\frac{2}{2-n}}(\exp_{i}(\frac{x}{\l}))}\right)^{\frac{n-2}{2}}=  
\left(
\frac
{1}
{1+\l^{2}G_{j}^{\frac{2}{2-n}}(a_{i})}
\right)^{\frac{n-2}{2}}
\left(
\frac
{1}
{1
+
\frac
{\l^{2}\left[G_{j}^{\frac{2}{2-n}}(\exp_{i}(\frac{x}{\l}))-G_{j}^{\frac{2}{2-n}}(a_{i})\right]}{1+\l^{2}G_{j}^{\frac{2}{2-n}}(a_{i})}}
\right)^{\frac{n-2}{2}} \\
& = 
\left(
\frac
{1}
{1+\l^{2}G_{j}^{\frac{2}{2-n}}(a_{i})}
\right)^{\frac{n-2}{2}}
\left(
\frac
{1}
{1
+
\frac
{\l\nabla G_{j}^{\frac{2}{2-n}}(a_{i})x+O(\vert x \vert^{2})}{1+\l^{2}G_{j}^{\frac{2}{2-n}}(a_{i})}}
\right)^{\frac{n-2}{2}} \\
& = 
\left(
\frac
{1}
{1+\l^{2}G_{j}^{\frac{2}{2-n}}(a_{i})}
\right)^{\frac{n-2}{2}}-\frac{(n-2)}{2}\l
\left(
\frac
{1}
{1+\l^{2}G_{j}^{\frac{2}{2-n}}(a_{i})}
\right)^{\frac{n}{2}}
\nabla G_{j}^{\frac{2}{2-n}}(a_{i})x
+
O\left(\frac{\vert x \vert^{2}}{(1+\l^{2}G_{j}^{\frac{2}{2-n}}(a_{i}))^{\frac{n}{2}}}\right).
\end{split}\end{equation}
Now, combining  \eqref{eq:41} and \eqref{eq:45}, we obtain
\begin{equation}\label{eq:47}
\begin{split}
\epsilon_{i,j}
= &
(1+o_{A_{i, j}, \l}(1))\frac{u_{j}(a_{i})c_3}{\left(1+\l^{2}G_{j}^{\frac{2}{2-n}}(a_{i})\right)^{\frac{n-2}{2}}}
+
o_{\varepsilon_{i, j}}(\varepsilon_{i,j}) \\
& +
\int_{\mathcal{A}^{c}\cap B^{\l}}
\left(\frac{1}{1+|x|^{2}}\right)^{\frac{n+2}{2}}
\left[\left(\frac{1}{1+\l^{2}G_{j}^{\frac{2}{2-n}}}\right)^{\frac{n-2}{2}}\right]\circ \exp_{i}(\frac{x}{\l})\;dx.
\end{split}\end{equation}
Next, using \eqref{eq:invarpro}, \eqref{eq:large1}, and Taylor expansion, we derive that
\begin{equation}\label{eq:48}
\begin{split}
\frac{u_{j}(a_{i})c_3}{\left(1+\l^{2}G_{j}^{\frac{2}{2-n}}(a_{i})\right)^{\frac{n-2}{2}}}
= &
c_3\left(1+o_{\varepsilon_{i, j}}(1)\right)u_{j}(a_{i})\frac{G_{j}(a_{i})}{\l ^{n-2}}\\
= &
c_3\left(1+o_{\varepsilon_{i, j}}(1)\right)\frac{G(a_{i},a_{j})}{\l^{n-2}}
=
\left(1+o_{\varepsilon_{i, j}}(1)\right)\varepsilon_{i,j}.
\end{split}
\end{equation}
Thus, combining \eqref{eq:equivalence0}, \eqref{eq:large1},  \eqref{eq:47}, and \eqref{eq:48}, we obtain
\begin{equation}\begin{split}
\epsilon_{i,j}
= &
\left(1+o_{\varepsilon_{i, j}}(1)\right)\varepsilon_{i,j} 
+
I_{\mathcal {A}^{c}},
\end{split}\end{equation}
where
\begin{equation}\begin{split}
I_{\mathcal{A}^c}=&
\int_{\mathcal{A}^{c}\cap B^{\l}(0)}
\left(\frac{1}{1+|x|^{2}}\right)^{\frac{n+2}{2}}
\left[\left(\frac{1}{1+\l^{2}G_{j}^{\frac{2}{2-n}}}\right)^{\frac{n-2}{2}}\right]\circ \exp_{i}(\frac{x}{\l})\;dx,
\end{split}\end{equation}
and
$$
\mathcal{A}^c:=\bar\R^n_+\setminus \mathcal{A}.
$$
Hence, to end the proof of ii) of point 1) and to prove iii) of point 1), we are going to show that \;$
I_{\mathcal{A}^c}$\; satisfies
\begin{equation}\label{eq:iaest}
I_{\mathcal{A}^c}=o_{\varepsilon_{i, j}}(\varepsilon_{i, j}).
\end{equation}
In order to do that, we first decompose \;$\mathcal{A}^{c}$\; into
\begin{equation}\label{eq:defb}
\begin{split}
\mathcal{B}
=
\{x\in \bar \R^n_+:\;\;
\gamma\l\sqrt{G_{j}^{\frac{2}{2-n}}(a_{i})}\leq \vert x \vert 
\leq 
\gamma^{-1}\l\sqrt{G_{j}^{\frac{2}{2-n}}(a_{i})}
\}
\end{split}\end{equation}
and
\begin{equation}\label{eq:deftildec}
\begin{split}
\mathcal{C}
=
\{x\in \bar\R^n_+:\;\;
\vert x \vert > \gamma^{-1}\l\sqrt{G_{j}^{\frac{2}{2-n}}(a_{i})}
\},
\end{split}\end{equation}
and  have
\begin{equation}\begin{split}
I_{\mathcal {A}^c}
=
I_{\mathcal{B}}
+
I_{\mathcal{C}},
\end{split}\end{equation}
where
\begin{equation}\label{eq:ib}
I_{\mathcal{B}}:=\int_{\mathcal{B}\cap B^{\l}(0)}
\left(\frac{1}{1+|x|^{2}}\right)^{\frac{n+2}{2}}
\left[\left(\frac{1}{1+\l^{2}G_{j}^{\frac{2}{2-n}}}\right)^{\frac{n-2}{2}}\right]\circ \exp_{i}(\frac{x}{\l})\;dx
\end{equation}
and
\begin{equation}\label{eq:ic}
I_{\mathcal{C}}:=\int_{\mathcal{C}\cap B^{\l}(0)}
\left(\frac{1}{1+|x|^{2}}\right)^{\frac{n+2}{2}}
\left[\left(\frac{1}{1+\l^{2}G_{j}^{\frac{2}{2-n}}}\right)^{\frac{n-2}{2}}\right]\circ \exp_{i}(\frac{x}{\l})\;dx.
\end{equation}
To prove \eqref{eq:iaest}, we are going to estimate separately $I_{\mathcal{B}}$ and $I_{\mathcal{C}}$. We start with $I_{\mathcal{B}}$. Using \eqref{eq:defb} and \eqref{eq:ib}, we have clearly that \;$I_{\mathcal{B}}$\; verifies the following estimate
\begin{equation}\label{eq:ibest1}
\begin{split}
I_{\mathcal{B}}
\leq 
\frac{C_{\gamma}}{\left(1+\l^{2}G_{j}^{\frac{2}{2-n}}(a_{i})\right)^{\frac{n+2}{2}}} \times
\int_{\mathcal{B}\cap B^{\l}}
\left[\left(\frac{1}{1+\l^{2}G_{j}^{\frac{2}{2-n}}}\right)^{\frac{n-2}{2}}\right]\circ \exp_{i}(\frac{x}{\l})\:dx,
\end{split}
\end{equation}
for some large positive constant $C_{\gamma}$ depending only on $\gamma$. Thus, rescaling and changing coordinates via $\exp_{j}\circ \exp_{i}^{-1}$ (if necessary), we have that \eqref{eq:ibest1} implies
\begin{equation}\label{eq:ibest2}
\begin{split}
I_{\mathcal{B}}
\leq & 
\hat C_{\gamma}\varepsilon_{i,j}^{\frac{n+2}{n-2}}
\int_{[\vert x \vert \leq \l\tilde C_{\gamma}d_{g}(a_{i},a_{j})]}
\left(\frac{1}{1+|x|^{2}}\right)^{\frac{n-2}{2}}dx
\leq
\bar C_{\gamma}\varepsilon_{i,j}^{\frac{n}{n-2}},
\end{split}\end{equation}
for some large positive constant $\hat C_{\gamma}$, $\tilde C_{\gamma}$, and $\bar C_{\gamma}$ which are depending only on $\gamma$.
Finally, we estimate $I_{\mathcal{C}}$. To do that, we fix $\gamma>0$ sufficiently small and use \eqref{eq:equivalence0}, \eqref{eq:large1}, and \eqref{eq:ic} to obtain
\begin{equation}\label{eq:icest1}
\begin{split}
I_{\mathcal{C}}
\leq 
\frac{C_{\gamma}}{\left(1+\l^{2}G_{j}^{\frac{2}{2-n}}(a_{i})\right)^{\frac{n-2}{2}}}
\int_{\mathcal{C}}
\left(\frac{1}{1+|x|^{2}}\right)^{\frac{n+2}{2}}dx
=
o_{\varepsilon_{i, j}}(\varepsilon_{i,j}),
\end{split}\end{equation}
for some large constant $C_{\gamma}$ depending only on $\gamma$. Hence \eqref{eq:ibest2} and \eqref{eq:icest1} implies \eqref{eq:iaest}, thereby ending the proof of point 1). On the other hand, we have clearly that point 2) follows from Lemma \ref{eq:bubbleinteractl}, Lemma \ref{eq:bsvolestl} and the fact that \;$B_g\varphi_{a, \l}=0$. Furthermore, the first equation of point 3) follows from Lemma \ref{eq:bubbleinteractl},  and ii) of point 1), while the second equation follows from the first equation and from the self-adjointness of \;$(L_g, B_g)$.
\end{pf}
\vspace{6pt}

\noindent
Now, using \eqref{eq:standardbubbles}, \eqref{eq:varphialb}, \eqref{eq:decompval}-\eqref{eq:walest}, we have the following interaction type estimate.
\begin{lem}\label{eq:alphabetainteract}
Assuming that $p\in \N^*$, $A:=(a_1, \cdots, a_p)\in (\partial M)^p$ and $\l\geq \frac{2}{\delta_0^{\frac{n}{2}}}$, then for every $i, j=1, \cdots, p$ with $i\neq j$, there holds
\begin{equation}\label{eq:betaint}
\begin{split} 
\int_M \varphi_{a_i, \l}^{\frac{n}{n-2}}\varphi_{a_j, \l}^{\frac{n}{n-2}}dV_g
=
O_{A_{i, j}, \l}(\varepsilon_{i,j}^{\frac{n}{n-2}}\log \varepsilon_{i,j}),
\end{split}
\end{equation} 
where $A_{i, j}=(a_i, a_j)$, $\varepsilon_{i, j}:=\varepsilon_{i, j}(A, \bar\l)$ with $\bar \l:=(\l, \cdots, \l)$, and for the meaning of  $O_{A_{i, j}, \l}(1)$ and $\epsilon_{i, j}(A, \bar\l)$, see respectively Section \ref{eq:notpre} and \eqref{eq:varepsilonij}.
\end{lem}
\begin{pf}
Using \eqref{eq:standardbubbles}, \eqref{eq:varphialb}, \eqref{eq:decompval}-\eqref{eq:walest} and setting \;$\varphi_i=\varphi_{a_i, \l }$\; and\; $\exp_i:=\exp_{a_i}^g$\; for \;$i=1, \cdots, p$ (for the meaning of $\exp_{a_i}^g$, see Section \ref{eq:notpre}), we have that for every \;$i=1, \cdots, p$, the following estimate holds
\begin{equation}\label{eq:betaint1}
\begin{split}
\varphi_{i} \leq C \left(\frac{\lambda}{1+\lambda^{2}d_{g}^{2}(a_{i},\cdot)}\right)^{\frac{n-2}{2}},
\end{split}\end{equation}
for some large positive constant independent of \;$a_i$\; and \;$\l$. Hence, using \eqref{eq:betaint1}, we have for \,$c$\; positive and small that the following estimate holds
\begin{equation}\label{eq:betaint2}
\begin{split}
\int_M \varphi_{i}^{\frac{n}{n-2}} \varphi_{j}^{\frac{n}{n-2}}dV_g
\leq &
C \underset{B_{0}^{\R^n_+}(c)}{\int}\left(\frac{\lambda }{1+\lambda ^2|x|^{2}}\right)^{\frac{n}{2}}\left(\frac{ \lambda  }{1 +\lambda ^{2} 
d_{g}^{2}(a_{j},\exp_{i}(x))}\right)^{\frac{n}{2}}dx\\
& +
C  \frac{1}{\lambda ^{\frac{n}{2}}}\underset{B_{0}^{\R^n_+}(c)}{\int}\left(\frac{ \lambda }{1+\lambda ^2|x|^{2}}\right)^{\frac{n}{2}}dx
+
O_{A_{i, j}, \l}\left(\frac{1}{\lambda ^{n}}\right) \\
= &
C \underset{B_{0}^{\R^n_+}(c\lambda)}{\int}\left(\frac{1}{1+r^{2}}\right)^{\frac{n}{2}}
\left(\frac{1}{1 
+
\lambda^{2}  d_{g}^{2}(a_{j},\exp_{i}(\frac{x}{\lambda}))}\right)^{\frac{n}{2}}dx\\
& +
C  \frac{1}{\lambda^{n}} 
\underset{B_{0}^{\R^n_+}(c\l)}{\int}\left(\frac{1}{1+|x|^{2}}\right)^{\frac{n}{2}}dx
+
O_{A_{i, j}, \l}\left(\frac{1}{\lambda^{n}}\right),
\end{split}\end{equation} 
for some large positive constant $C$ independent of $A_{i, j}$ and $\l$. So, appealing to \eqref{eq:betaint2}, we infer that
\begin{equation}\label{eq:betaint3}
\begin{split}
\int_M \varphi_{i}^{\frac{n}{n-2}}  \varphi_{j}^{\frac{n}{n-2}}dV_g 
\leq &
C\underset{B_{0}^{\R^n_+}(c\l)}{\int}\left(\frac{1}{1+|x|^{2}}\right)^{\frac{n}{2}}\left(\frac{1}{1 +\lambda^{2}  
 d_{g}^{2}(a_{j},\exp_{i}(\frac{x}{\lambda}))}\right)^{\frac{n}{2}}dx
+
O_{A_{i, j}, \l}(\frac{\log \lambda}{\lambda^{n}}).
\end{split}\end{equation} 
Thus \eqref{eq:betaint} follows from \eqref{eq:betaint3} if 
\begin{equation}\begin{split}
d_g( a _{i}, a _{j})\geq 3c.
\end{split}\end{equation} 
Hence to complete the proof of the lemma it remains to treat the case $d_g( a _{i}, a _{j})<3c$. To do that, we set 
\begin{equation*}
\mathcal{B}=\{x\in \bar \R^n_+:\;\;\;\;\frac{1}{2}d_{g}( a _{i}, a _{j})\leq \vert  \frac{x}{\lambda }\vert \leq 2 d_{g}( a _{i}, a _{j})\},
\end{equation*}
and use \eqref{eq:betaint3} combined with the triangle inequality to get for $c>0$ sufficiently small that the following estimate holds
\begin{equation}\begin{split}
\int_M \varphi_{i}^{\frac{n}{n-2}}&\varphi_{j}^{\frac{n}{n-2}}dV_g
\leq 
C \underset{\mathcal{B}}{\int}
\left(\frac{1}{1+|x|^{2}}\right)^{\frac{n}{2}}\left(\frac{1}{1 +\lambda^{2} d_{g}^{2}( a _{j}, \exp_{i}( \frac{x}{\lambda }))}\right)^{\frac{n}{2}}dx +
O_{A_{i, j}, \l}(\varepsilon_{i,j}^{\frac{n}{n-2}}\log \varepsilon_{i,j})\\
\leq &
C
\left(\frac{1}{1+\vert\lambda d_{g}( a _{i}, a _{j})\vert^{2}}\right)^{\frac{n}{2}}
\int_{\{\vert  \frac{x}{\lambda }\vert \leq 4 d_g( a _{i}, a _{j})\}}
\left(\frac{1}{1+|x|^{2}}\right)^{\frac{n}{2}}dx 
+
O_{A_{i, j}, \l}(\varepsilon_{i,j}^{\frac{n}{n-2}}\log \varepsilon_{i,j})\\
= &
O_{A_{i, j}, \l}(\varepsilon_{i,j}^{\frac{n}{n-2}}\log \varepsilon_{i,j}),
\end{split}
\end{equation} 
where $C$ is a large positive constant independent of $A_{i, j}$ and $\l$, thereby completing the proof of the lemma.
\end{pf}
\vspace{6pt}

\noindent
To end this section,  we extend the definition of $\varphi_{a, \l}$ to interior concentration points by just using the classical compactification-type bubbles of Schoen\cite{schoen}. Precisely, for $a\in M$ and $0<2\epsilon\leq \delta\leq \delta_0$, we set 
\begin{equation}\label{eq:varphiali} 
v_{a, \epsilon, \delta}(\cdot):=\chi_{\d}(\cdot)\delta_{a, \epsilon}(\cdot)+\left(1-\chi_{\delta}(\cdot)\right)\epsilon^{\frac{n-2}{2}}G_a(a, \cdot),
\end{equation}
$\chi_{\d}$ is as in \eqref{eq:cutoff}, $\delta_{a, \epsilon}$ is defined as in \eqref{eq:standardbubbles}, $G_{a}(a, \cdot)$ is defined as in \eqref{eq:greena}, and for \;$\l\geq \frac{2}{\d_0^{\frac{n}{2}}}$, we define
\begin{equation}\label{eq:varphialint}
\varphi_{a, \l}:=u_av_{a,\l, (\frac{1}{\l})^{\frac{2}{n}}},
\end{equation}
where \;$u_a$\; is as \eqref{eq:detga} and \;$\d_0$\; is as in Lemma \ref{eq:brenchenenergy}.
\begin{rem}\label{eq:notsharpintbubble}
We would like to point out the fact that for bubbles which are centered at interior points, we do not need sharp ones, but just bubbles which can be used to describe diverging PS sequences of\; $\mathcal{E}_g$. This is due to the fact that our scheme of the algebraic topological argument of Bahri-Coron\cite{bc} relies only on bubbles centered at boundary points, see Section \ref{eq:algtop}.
\end{rem}

\section{Mapping \;$B_p(\partial M)$\; into appropriate sublevels of\;$\mathcal{E}_g$}\label{eq:mappingbary}
In this section, we map $B_p(\partial M)$ into some appropriate sublevels of the Euler-Lagrange functional $\mathcal{E}_g$ via the Brendle-Chen\cite{bs}'s bubbles. Precisely, we are going to derive sharp energy estimates for convex combination of the bubbles $\varphi_{a, \l}$ given by \eqref{eq:varphialb} so that we can confirm in the next section that they are algebraic topological ones by carrying a suitable scheme of the  barycenter technique of Bahri-Coron\cite{bc}. In order to do that, we first make the following definition. For \;$p\in \N^*$, $\sigma:=\sum_{i=1}^p\alpha_i\d_{a_i}\in B_p(\partial M)$\; and \;$\l\geq \frac{2}{\delta_0^{\frac{n}{2}}}$ with \;$\d_0$\; given by Lemma \ref{eq:brenchenenergy}, we define \;$f_p(\l): \;B_p(\partial M)\longrightarrow W^{1, 2}_+(\ov M)$ as follows 
\begin{equation}\label{eq:deffp}
f_p(\l)(\sigma):=\sum_{i=1}^p\alpha_i\varphi_{a_i, \l}.
\end{equation}
Now, we start the goal of this section with the following proposition which provides the first step to apply our scheme of the algebraic topological argument of Bahri-Coron\cite{bc}.
\begin{pro}\label{eq:baryest}
 There exists a large constant \;$C_0>0$, $\nu_0>1$\; and \;$0<\varepsilon_0\leq \delta_0$\;such that for every \;$p\in \N^*$ and every \;$0<\varepsilon\leq \varepsilon_0$, there exists \;$\l_p:=\l_p(\varepsilon):=\l_p(\nu_0, \varepsilon)\geq \frac{2}{\d_0^{\frac{n}{2}}}$\; such that for every $\l\geq \l_p$ and for every $\sigma=\sum_{i=1}^p\alpha_i\delta_{a_i}\in B_p(\partial M)$, we have\\\\
 1) If there exist \;$i_0\neq j_0$\; such that \;$\frac{\alpha_{i_0}}{\alpha_{j_0}}>\nu_0$\; or if \;$\sum_{i\neq j}\varepsilon_{i, j}> \varepsilon$, then
 $$
\mathcal{E}_g(f_p(\l))(\sigma))\leq p^{\frac{2}{n}}\mathcal{Y}(S^n_+),
$$ 
where \;$\mathcal{Y}(S^n_+)$\; is defined by \eqref{eq:yamabehsphere} and \;$\varepsilon_{i, j}:=\varepsilon_{i, j}(A, \bar \l)$ with $\bar \l:=(\l, \cdots, \l)$ and for the definition of \;$\varepsilon_{i, j}(A, \bar \l)$, see \eqref{eq:varepsilonij}.\\
2) If for every $i\neq j$ we have \;$\frac{\alpha_{i}}{\alpha_j}\leq\nu_0$\; and if \;$\sum_{i\neq j}\varepsilon_{i, j}\leq \varepsilon$, then
$$
\mathcal{E}_g(f_p(\l))(\sigma))\leq p^{\frac{2}{n}}\mathcal{Y}(S^n_+)\left(1+\frac{C_0}{\l^{n-2}}-c_g\frac{(p-1)}{\l^{n-2}}\right),
$$
where \;$c_g$\; is is defined by \eqref{eq:defcg}.
\end{pro}
\vspace{6pt}

\noindent
Proposition \ref{eq:baryest} is derived from the following technical Lemma.
\begin{lem}\label{eq:baryestaux}
We have that the following holds:\\
1)
For every \;$\epsilon$\; positive and small and for every \;$p\in \N^*$, there exists\; $\lambda_{p}:=\l_p(\epsilon)\geq \frac{2}{\d_0^{\frac{n}{2}}}$ such that for every \;$\l\geq \l_p$ and for every $\sigma:=\sum_{i=1}^p\alpha_i\d_{a_i}\in B_p(\partial M)$, we have
\begin{equation*}
\begin{split}
\sum_{i\neq j}\epsilon_{i,j}>\epsilon
\end{split}
\end{equation*}
implies
\begin{equation*}
\begin{split}
\mathcal{E}_{g}(f_p(\l)(\sigma))< p^{\frac{2}{n}}\mathcal{Y}(S^n_+),
\end{split}
\end{equation*}
where \;$\epsilon_{i, j}:=\epsilon_{i, j}(A, \l)$ is defined by \eqref{eq:intweight}.\\
2) For every \;$\nu>1$, for every \;$\epsilon>0$\;and small, and for every \;$p\in \N^*$, there exists \;$\lambda_{p}:=\lambda_{p}(\epsilon,\nu)\geq \frac{2}{\d_0^{\frac{n}{2}}}$\; such that for every $\l\geq \l_p$ and for every \;$\sigma:=\sum_{i=1}^p\alpha_i\d_{a_i}\in B_p(\partial M)$, we have
\begin{equation*}\begin{split}
\exists\;{i_0\neq j_0}\;\;\text{such that}\;\;\;\frac{\alpha_{i_0}}{\alpha_{j_0}} >\nu
\; \;\;\;\text{and}\; \;\;\;\sum_{i\neq j}\epsilon_{i,j}\leq\epsilon
\end{split}
\end{equation*}
imply
\begin{equation*}
\begin{split}
\mathcal{E}_{g}(f_p(\l)(\sigma))< p^{\frac{2}{n}}\mathcal{Y}(S^n_+).
\end{split}
\end{equation*}
3) There exists \;$C_0>0$, \;$\nu_0>1$, \;$\l_0\geq \frac{2}{\d_0^{\frac{n}{2}}}$\; and \;$0<\epsilon_0\leq\delta_0$\; such that for every \;$1<\nu\leq \nu_0$, for every \;$0<\epsilon\leq\epsilon_0$, for every \;$p\in \N^*$, for every \;$\l\geq \l_0$, and for every $\sigma:=\sum_{i=1}^p\alpha_i\d_{a_i}\in B_p(\partial M)$, we have
\begin{equation}\begin{split}
\frac{\alpha_{i}}{\alpha_{j}} \leq\nu\;\;\;\forall i, j,\;\;
\; \text{and}\;\;\; \sum_{i\neq j}\epsilon_{i,j}\leq\epsilon
\end{split}\end{equation}
imply
\begin{equation}\begin{split} 
\mathcal{E}_{g}(f_p(\l)(\s))
\leq &
p^{\frac{2}{n}}\mathcal{Y}(S^n_+)
\left(
1+\frac{C_0}{\lambda^{n-2}}
-
c_g\frac{(p-1)}{\lambda^{n-2}}
\right).
\end{split}\end{equation} 
\end{lem}
\begin{pf}
First of all, we set 
\begin{equation}\label{eq:numden}
 \mathcal{N}_g(u):=\langle L_gu,u\rangle+\langle B_gu, u\rangle, \;\;\mathcal{D}_g(u):=\left(\int_M u^{\frac{2n}{n-2}}dV_g\right)^{\frac{n-2}{n}}, \;\;\;\;u\in W^{1, 2}_+(\ov M),
\end{equation}
and use \eqref{eq:Yamabefunctional} to have
\begin{equation}\label{eq:eulerlagrange}
\mathcal{E}_g(u)=\frac{\mathcal{N}_g(u)}{\mathcal{D}_g(u)},\;\;\;\;u\in W^{1, 2}_+(\ov M).
\end{equation}
Furthermore, for $p\in \N^*$, $\sigma:=\sum_{i=1}\alpha_i\d_{a_i}\in B_p(\partial M)$ and $\l\geq \frac{2}{\d_0^{\frac{n}{2}}}$, we set (as in the proof of Lemma \ref{eq:interactestl})
\begin{equation}\label{eq:varphii1}
 \varphi_i=\varphi_{a_i, \l}, \;\;i=1, \cdots, p.
\end{equation}
Now, we start with the proof of point 1). To do so, we first use Lemma \ref{eq:interactestl}, \eqref{eq:deffp}, \eqref{eq:numden}, \eqref{eq:varphii1}, and H\"older's inequality to estimate $\mathcal{N}_g(f_p(\l)(\s))$ as follows
\begin{equation}\label{eq:ngest1}
\begin{split} 
\mathcal{N}_{g}(f_p(\l)(\s))
=&
c_{0}(1+o_{A, \bar\alpha, \l}(1))
\int_M \left(\sum_{i=1}^p\alpha_{i}\varphi_{i}^{\frac{n+2}{n-2}}\right)\left(\sum_{j=1}^p\alpha_{j}\varphi_{j}\right)dV_g\\
= &
c_{0}(1+o_{A, \bar\alpha, \l}(1))
\int_M \left(\frac{\sum_{i=1}^p\alpha_{i}\varphi_{i}^{\frac{n+2}{n-2}}}{ \sum_{j=1}^p\alpha_{j}\varphi_{j}}\right)\left(\sum_{j=1}^p\alpha_{j}\varphi_{j}\right)^{2}dV_g \\
\leq &
c_{0}(1+o_{A, \bar\alpha, \l}(1))
D_{g}(u)\Vert\frac{\sum_{i=1}^p\alpha_{i}\varphi_{i}^{\frac{n+2}{n-2}}}{\sum_{j=1}^p\alpha_{j}\varphi_{j}}\Vert_{L^{\frac{n}{2}}},
\end{split}
\end{equation} 
where \;$A:=(a_1, \cdots, a_p)$, $\bar\alpha:=(\alpha_1, \cdots, \alpha_p)$\; and for the meaning of \;$o_{A, \bar\alpha, \l}(1)$, see Section \ref{eq:notpre}.
Thus, using the convexity of the map $x\longrightarrow x^{\beta}$ with $\beta>1$, we derive that \eqref{eq:ngest1} implies  
\begin{equation}\label{eq:ngest2}
\begin{split} 
\mathcal{N}_{g}(f_p(\l)(\s))
\leq &
c_{0}(1+o_{A, \bar\alpha, \l}(1))
\mathcal{D}_{g}(u)\left(\int_M\left(\sum_{i=1}^p\frac{\alpha_{i}\varphi_{i}}{\sum_{j=1}^p\alpha_{j}\varphi_{j}}\varphi_{i}^{\frac{4}{n-2}}\right)^{\frac{n}{2}}\right)^{\frac{2}{n}}dV_g \\
\leq &
c_{0}(1+o_{A, \bar\alpha, \l}(1))
\mathcal{D}_{g}(u)\left(\sum_{i=1}^p\int_M\frac{\alpha_{i}\varphi_{i}}{\sum_{j=1}^p\alpha_{j}\varphi_{j}}\varphi_{i}^{\frac{2n}{n-2}}\right)^{\frac{2}{n}}dV_g.
\end{split}\end{equation} 
Hence, clearly Lemma \ref{eq:bsvolestl}, \eqref{eq:eulerlagrange} and \eqref{eq:ngest2} imply for any pair \;$i\neq j$ ($i, j=1, \cdots, p$)
\begin{equation}\label{eq:ngest3}
\begin{split} 
\mathcal{E}_{g}(f_p(\l)(\s))
\leq & 
c_{0}(1+o_{A, \bar\alpha, \l}(1))
\left(
c_{1}(p-1)
+
\int_M \frac{\alpha_{i}\varphi_{i}}{\alpha_{i}\varphi_{i}+\alpha_{j}\varphi_{j}}\varphi_{i}^{\frac{2n}{n-2}}dV_g
\right)^{\frac{2}{n}} \\
\leq &
c_{0}(1+o_{A, \bar\alpha,  \l}(1))
\left(
c_{1}p
-
\int_M \frac{\alpha_{j}\varphi_{j}}{\alpha_{i}\varphi_{i}+\alpha_{j}\varphi_{j}}\varphi_{i}^{\frac{2n}{n-2}}
dV_g\right)^{\frac{2}{n}}, \\
\end{split}\end{equation} 
and we may assume \;$\alpha_{i}\leq\alpha_{j}$\; by symmetry. Now, we are going to estimate from below the quantity $\int_M \frac{\alpha_{j}\varphi_{j}}{\alpha_{i}\varphi_{i}+\alpha_{j}\varphi_{j}}\varphi_{i}^{\frac{2n}{n-2}}dV_g
$. In order to do that, for $\gamma>0$, we set
\begin{equation}\label{eq:defmathaij}
\begin{split} 
\mathcal{A}_{i,j}=
\{x\in M:\;\;\;
\varphi_{i}(x)\geq \gamma (\frac{\alpha_{i}}{\alpha_{j}}\varphi_{i}(x)+\varphi_{j}(x))
\},
\end{split}\end{equation} 
and use \eqref{eq:defmathaij} to have
\begin{equation}\label{eq:ngest4}
\begin{split} 
\int_M \frac{\alpha_{j}\varphi_{j}\varphi_{i}^{\frac{2n}{n-2}}}{\alpha_{i}\varphi_{i}+\alpha_{j}\varphi_{j}}dV_g
\geq &
\int_{\mathcal{A}_{i,j}} \frac{\varphi_{j}}{\frac{\alpha_{i}}{\alpha_{j}}\varphi_{i}+\varphi_{j}}\varphi_{i}^{\frac{2n}{n-2}}dV_g
\geq 
\gamma 
\int_{\mathcal{A}_{i,j}}
\varphi_{i}^{\frac{n+2}{n-2}}\varphi_{j}dV_g \\
= &
\gamma 
\left(
\int_M \varphi_{i}^{\frac{n+2}{n-2}}\varphi_{j}dV_g-\int_{\mathcal{A}_{i,j}^{c}}\varphi_{i}^{\frac{n+2}{n-2}}\varphi_{j}dV_g
\right) \\
\geq &
\gamma 
\left(
\int_M \varphi_{i}^{\frac{n+2}{n-2}}\varphi_{j}dV_g
-
\gamma ^{\frac{4}{n-2}}
\int_{\mathcal{A}_{i,j}^{c}}\left(\frac{\alpha_{i}}{\alpha_{j}}\varphi_{i}+\varphi_{j}\right)^{\frac{4}{n-2}}\varphi_{i}\varphi_{j}dV_g
\right),
\end{split}\end{equation} 
where $\mathcal{A}_{i,j}^c:=M\setminus \mathcal{A}_{i, j}$.
Next, since $\frac{\alpha_{i}}{\alpha_{j}}\leq 1$, then appealing to \eqref{eq:ngest4}, we infer that the following estimate holds 
\begin{equation}\label{eq:ngest5}
\begin{split} 
\int_M \frac{\alpha_{j}\varphi_{j}\varphi_{i}^{\frac{2n}{n-2}}}{\alpha_{i}\varphi_{i}+\alpha_{j}\varphi_{j}}dV_g
\geq &
\gamma 
\left(
\int_M \varphi_{i}^{\frac{n+2}{n-2}}\varphi_{j}dV_g
-
C\gamma ^{\frac{4}{n-2}}
\int_M(\varphi_{i}^{\frac{4}{n-2}}+\varphi_{j}^{\frac{4}{n-2}})\varphi_{i}\varphi_{j}dV_g
\right),
\end{split}\end{equation} 
for some large positive constant \;$C$\; independent of $A$, $\l$ and $\gamma$. Thus, ii) of point 1) of Lemma \ref{eq:interactestl} and  \eqref{eq:ngest5} imply that for \;$\gamma>0$\; sufficiently small, there holds
\begin{equation}\label{eq:ngest6}
\begin{split} 
\int_M \frac{\alpha_{j}\varphi_{j}}{\alpha_{i}\varphi_{i}+\alpha_{j}\varphi_{j}}\varphi_{i}^{\frac{2n}{n-2}}dV_g
\geq &
\frac{\gamma }{2}
\int_M \varphi_{i}^{\frac{n+2}{n-2}}\varphi_{j}dV_g.
\end{split}\end{equation} 
Hence, combining \eqref{eq:ngest3} and \eqref{eq:ngest6}, we conclude that for any pair \;$i\neq j$, the following estimate holds
\begin{equation}\label{eq:ngest7}
\begin{split} 
\mathcal{E}_{g}(f_p(\l)(\s))
\leq & 
\left(1+o_{A, \bar\alpha, \l}(1)\right)\mathcal{Y}(S^n_+)
\left(
p
-
\frac{\gamma }{2c_1}
\int_M \varphi_{i}^{\frac{n+2}{n-2}}\varphi_{j}dV_g
\right)^{\frac{2}{n}}.
\end{split}\end{equation} 
Clearly \eqref{eq:ngest7} implies, that we always have 
\begin{equation}\begin{split} 
\mathcal{E}_g(f_p(\l)(\s))\leq \left(1+o_{A, \bar\alpha, \l}(1)\right)p^{\frac{2}{n}}\mathcal{Y}(S^n_+)
\end{split}\end{equation} 
and in case $\sum_{i\neq j}\epsilon_{i,j}\geq \epsilon$ 
\begin{equation}\begin{split} 
\mathcal{E}_{g}(f_p(\l)(\s))
\leq & 
\left(1+o_{A, \bar\alpha, \l}(1)\right)p^{\frac{2}{n}}\mathcal{Y}(S^n_+)
\left(
1
-
\frac{\gamma \epsilon}{2pc_1}
\right)^{\frac{2}{n}}.
\end{split}\end{equation} 
thereby ending the proof of point 1). Now, we are going to treat the second case. Hence, we may assume
\begin{equation}\begin{split} 
\sum_{i\neq j}\epsilon_{i,j}\ll 1
\end{split}\end{equation} 
and thus according to Lemma \ref{eq:interactestl}
\begin{equation}\label{eq:varepijep}
\begin{split}
\epsilon_{i,j}=(1+o_{\varepsilon_{i, j}}(1))\varepsilon_{i,j}
\;\text{ and }\;
\lambda d_{g}(a_{i},a_{j})\gg 1,
\end{split}
\end{equation}
and for the meaning of \;$o_{\epsilon_{i, j}}(1)$, see Section \ref{eq:notpre}. We then use Lemma \ref{eq:interactestl},  \eqref{eq:numden}, and \eqref{eq:varepijep} to have 
\begin{equation}\begin{split} \label{eq:ngest8}
\mathcal{N}_{g}(f_p(\l)(\s))
= &
\sum_{i=1}^p\sum_{j=1}^p\alpha_{i}\alpha_{j}\left(\langle L_{g}\varphi_{i}, \varphi_{j}\rangle+\langle B_g \varphi_i, \varphi_j\rangle\right) \\
= &
\sum_{i=1}^p\alpha_{i}^{2}\left(\langle L_{g}\varphi_{i}, \varphi_{i}\rangle+\langle B_g \varphi_i, \varphi_i \rangle \right)
+
\sum_{i\neq j}\alpha_{i}\alpha_{j}\left(\langle L_{g}\varphi_{i}, \varphi_{j}\rangle+\langle B_g\varphi_i, \varphi_j\rangle \right) \\
= &
\sum_{i}\alpha_{i}^{2}\mathcal{E}_{g}(\varphi_{i})\mathcal{D}_{g}(\varphi_{i})
+
c_{0}(1+o_{\sum_{i,\neq j}\varepsilon_{i, j}}(1))\sum_{i\neq j}\alpha_{i}\alpha_{j}\varepsilon_{i,j}
\end{split}\end{equation} 
and (for the meaning of \;$o_{\sum_{i\neq j}\epsilon_{i, j}}(1)$, see Section \ref{eq:notpre})
\begin{equation}\label{eq:dgest1}
\begin{split} 
\mathcal{D}^{\frac{n}{n-2}}_{g}(f_p(\l)(\s))
= &
\int_M \left(\sum_{i=1}^p\alpha_{i}\varphi_{i}\right)^{\frac{2n}{n-2}}dV_g
=
\sum_{i=1}^p\alpha_{i}\int_M \left(\sum_{j=1}^p\alpha_{j}\varphi_{j}\right)^{\frac{n+2}{n-2}}\varphi_{i}dV_g \\
= &
\sum_{i=1}^p\alpha_{i}\int_M \left(\alpha_{i}\varphi_{i}+\sum_{j=1, \;j\neq i}^{p}\alpha_{j}\varphi_{j}\right)^{\frac{n+2}{n-2}}\varphi_{i}dV_g.
\end{split}\end{equation} 
To proceed further, we set 
$
\mathcal{A}_{i}=\{x\in M:\ \;\;\;\alpha_{i}\varphi_{i}(x)>\sum_{j=1, \;j\neq i}^{p}\alpha_{j}\varphi_{j}(x)\},
$
and use Taylor expansion to obtain
\begin{equation}\label{eq:dgest2}
\begin{split} 
\mathcal{D}^{\frac{n}{n-2}}_{g}(f_p(\l)(\s))
= &
\sum_{i=1}^p\alpha_{i}^{\frac{2n}{n-2}}\int_{\mathcal{A}_{i}} \varphi_{i}^{\frac{2n}{n-2}}dV_g
+
\frac{n+2}{n-2}\sum_{i\neq j}\alpha_{i}^{\frac{n+2}{n-2}}\alpha_{j}\int_{\mathcal{A}_{i}} \varphi_{i}^{\frac{n+2}{n-2}}\varphi_{j}dV_g\\
+ &
\sum_{i=1}^p\alpha_{i}\int_{\mathcal{A}_{i}^{c}}\left(\sum_{ j=1, \;j\neq i}^{p}\alpha_{j}\varphi_{j}\right)^{\frac{n+2}{n-2}}\varphi_{i}dV_g +
O_{A, \bar\alpha, \l}
\left(\sum_{i\neq j}\alpha_{i}^{\frac{n}{n-2}}\alpha_{j}^{\frac{n}{n-2}}\int_M \varphi_{i}^{\frac{n}{n-2}}\varphi_{j}^{\frac{n}{n-2}}dV_g\right),
\end{split}\end{equation} 
where $\mathcal{A}_i^c=M\setminus \mathcal{A}_i$, \;$A:=(a_1, \cdots, a_p)$, $\bar\alpha:=(\alpha_1, \cdots, \alpha_p)$, $O_{A, \bar\alpha, \l}(1)$ is defined as in Section \ref{eq:notpre}, and we made use of $n\geq 4$ and the algebraic relation
\begin{equation}\begin{split}
(a+b)^{\frac{n}{n-2}}\leq C_n(a^{\frac{n}{n-2}}+b^{\frac{n}{n-2}})
\end{split}\end{equation}
for $a,b\geq 0$ and $C_n$ a positive constant depending only on $n$. Moreover since
\begin{equation}\begin{split} 
(a+b)^{\frac{n+2}{n-2}}\geq a^{\frac{n+2}{n-2}}+b^{\frac{n+2}{n-2}},
\end{split}\end{equation} 
then \eqref{eq:dgest2} implies
\begin{equation}\label{eq:dgest3}
\begin{split} 
\mathcal{D}^{\frac{n}{n-2}}_{g}(f_p(\l)(\s))
= &
\sum_{i=1}^p\alpha_{i}^{\frac{2n}{n-2}}\int_M \varphi_{i}^{\frac{2n}{n-2}}dV_g 
+
\frac{n+2}{n-2}\sum_{i\neq j}\alpha_{i}^{\frac{n+2}{n-2}}\alpha_{j}\int_M \varphi_{i}^{\frac{n+2}{n-2}}\varphi_{j}dV_g \\
& + 
\sum_{i\neq j}\alpha_{i}\alpha_{j}^{\frac{n+2}{n-2}}\int_{A_{i}^{c}} \varphi_{j}^{\frac{n+2}{n-2}}\varphi_{i}dV_g \\
& +
O_{A, \l}
\left(\sum_{i\neq j}\alpha_{i}^{\frac{n}{n-2}}\alpha_{j}^{\frac{n}{n-2}}\int_M \varphi_{i}^{\frac{n}{n-2}}\varphi_{j}^{\frac{n}{n-2}}dV_g\right) \\
= &
\sum_{i=1}^p\alpha_{i}^{\frac{2n}{n-2}}\int_M \varphi_{i}^{\frac{2n}{n-2}}dV_g 
+
\frac{2n}{n-2}\sum_{i\neq j}\alpha_{i}^{\frac{n+2}{n-2}}\alpha_{j} \int_M \varphi_{i}^{\frac{n+2}{n-2}}\varphi_{j}dV_g
 \\
& +
O_{A, \bar\alpha, \l}
\left(\sum_{i\neq j}\alpha_{i}^{\frac{n}{n-2}}\alpha_{j}^{\frac{n}{n-2}}\int_M \varphi_{i}^{\frac{n}{n-2}}\varphi_{j}^{\frac{n}{n-2}}dV_g\right).
\end{split}\end{equation} 
So, using Lemma \ref{eq:alphabetainteract} and \eqref{eq:varepijep}, we have that \eqref{eq:dgest3} implies
\begin{equation}\label{eq:dgest4}
\begin{split} 
\mathcal{D}^{\frac{n}{n-2}}_{g}(f_p(\l)(\s))
= &
\sum_{i=1}^p\alpha_{i}^{\frac{2n}{n-2}}\int_M \varphi_{i}^{\frac{2n}{n-2}}dV_g  \\
& +
\frac{2n}{n-2}\left(1+o_{\sum_{i\neq j}\varepsilon_{i, j}}(1)\right)\sum_{i\neq j}\alpha_{i}^{\frac{n+2}{n-2}}\alpha_{j} \varepsilon_{i,j}
+o_{\sum_{i\neq j}\varepsilon_{i, j}}\left(\sum_{i\neq j}\alpha_{i}^{\frac{n}{n-2}}\alpha_{j}^{\frac{n}{n-2}}\varepsilon_{i,j}\right),
\end{split}
\end{equation} 
Thus, using Young's inequality and the symmetry of $\varepsilon_{i,j}$, we infer from \eqref{eq:dgest4} that the following estimate holds
\begin{equation}\label{eq:dgest5}
\begin{split} 
\mathcal{D}^{\frac{n}{n-2}}_{g}(f_p(\l)(\s)))
= &
\sum_{i=1}^p\alpha_{i}^{\frac{2n}{n-2}}\int_M \varphi_{i}^{\frac{2n}{n-2}} dV_g 
+
\frac{2n}{n-2}(1+o_{\sum_{i\neq j}\varepsilon_{i, j}}(1))\sum_{i\neq j}\alpha_{i}^{\frac{n+2}{n-2}}\alpha_{j} \varepsilon_{i,j}.
\end{split}\end{equation} 
Hence, using again Young's inequality, Taylor expansion, and Lemma \ref{eq:bsvolestl}, we have that \eqref{eq:dgest5} gives
\begin{equation}\label{eq:dgest6}
\begin{split} 
\mathcal{D}_{g}(f_p(\l)(\s))
=& 
\left(\sum_{i=1}^p\alpha_{i}^{\frac{2n}{n-2}}\int_M \varphi_{i}^{\frac{2n}{n-2}}dV_g\right)^{\frac{n-2}{n}}\times
\left(
1
+
\frac
{
\frac{2n}{n-2}\left(1+o_{\sum_{i\neq j}\varepsilon_{i, j}}(1)\right)\sum_{i\neq j}\alpha_{i}^{\frac{n+2}{n-2}}\alpha_{j} \varepsilon_{i,j}
}
{
\sum_{i=1}^p\alpha_{i}^{\frac{2n}{n-2}}\int_M \varphi_{i}^{\frac{2n}{n-2}}dV_g
}
\right)^{\frac{n-2}{n}} 
 \\
= &
\left(\sum_{i=1}^p\alpha_{i}^{\frac{2n}{n-2}}\int_M \varphi_{i}^{\frac{2n}{n-2}}dV_g\right)^{\frac{n-2}{n}}
+
2
\frac
{
\left(1+o_{\sum_{i\neq j}\varepsilon_{i, j}}(1)\right)\sum_{i\neq j}\alpha_{i}^{\frac{n+2}{n-2}}\alpha_{j} \varepsilon_{i,j}
}
{
\left(\sum_{i=1}^p\alpha_{i}^{\frac{2n}{n-2}}\int_M \varphi_{i}^{\frac{2n}{n-2}}dV_g\right)^{\frac{2}{n}}
}
\\
= &
\left(\sum_{i=1}^p\alpha_{i}^{\frac{2n}{n-2}}\mathcal{D}_{g}^{\frac{n}{n-2}}(\varphi_{i})\right)^{\frac{n-2}{n}}
+
2c_{1}^{-\frac{2}{n}}
\frac
{
\left(1+o_{\sum_{i\neq j}\varepsilon_{i, j}}(1)\right)\sum_{i\neq j}\alpha_{i}^{\frac{n+2}{n-2}}\alpha_{j} \varepsilon_{i,j}
}
{
\left(\sum_{i=1}^p\alpha_{i}^{\frac{2n}{n-2}}\right)^{\frac{2}{n}}
}.
\end{split}\end{equation} 
Now, combining \eqref{eq:ngest8} and \eqref{eq:dgest6}, and using again Taylor expansion,  we obtain 
\begin{equation}\label{eq:ngest9}
\begin{split} 
\mathcal{E}_{g}(f_p(\l)(\s))
\leq &
\frac
{
\sum_{i=1}^p\alpha_{i}^{2}\mathcal{E}_{g}(\varphi_{i})\mathcal{D}_{g}(\varphi_{i})
+
c_{0}\left(1+o_{\sum_{i\neq j}\varepsilon_{i, j}}(1)\right)\sum_{i\neq j}\alpha_{i}\alpha_{j}\varepsilon_{i,j}
}
{
\left(\sum_{i=1}^p\alpha_{i}^{\frac{2n}{n-2}}\mathcal{D}_{g}^{\frac{n}{n-2}}(\varphi_{i})\right)^{\frac{n-2}{n}}
+
2c_{1}^{-\frac{2}{n}}
\frac
{
\left(1+o_{\sum_{i\neq j}\varepsilon_{i, j}}(1)\right)\sum_{i\neq j}\alpha_{i}^{\frac{n+2}{n-2}}\alpha_{j} \varepsilon_{i,j}
}
{
\left(\sum_{i=1}^p\alpha_{i}^{\frac{2n}{n-2}}\right)^{\frac{2}{n}}
}
} \\
= &
\frac
{
\sum_{i=1}^p\mathcal{E}_{g}(\varphi_{i})\alpha_{i}^{2}\mathcal{D}_{g}(\varphi_{i})
}
{
\left(\sum_{i=1}^p\alpha_{i}^{\frac{2n}{n-2}}\mathcal{D}_{g}^{\frac{n}{n-2}}(\varphi_{i})\right)^{\frac{n-2}{n}}
} +
\frac{c_{0}}{c_{1}^{\frac{n-2}{n}}}\frac{\left(1+o_{\sum_{i\neq j}\varepsilon_{i, j}}(1)\right)\sum_{i\neq j}\alpha_{i}\alpha_{j}\varepsilon_{i,j}}{\left(\sum_{i=1}^p\alpha_{i}^{\frac{2n}{n-2}}\right)^{\frac{n-2}{n}}} \\
& -
\frac{2c_{0}}{c_{1}^{\frac{n-2}{n}}}\frac{\left(1+o_{\sum_{i\neq j}\varepsilon_{i, j}}(1)\right)\left(\sum_{i=1}^p\alpha_{i}^{2}\right)\left(\sum_{i\neq j}\alpha_{i}^{\frac{n+2}{n-2}}\alpha_{j}\varepsilon_{i,j}\right)}{(\sum_{i=1}^p\alpha_{i}^{\frac{2n}{n-2}})^{\frac{2n-2}{n}}}.
\end{split}\end{equation} 
Hence, using \eqref{eq:relationcy} and rearranging the terms in \eqref{eq:ngest9}, we get
\begin{equation}\label{eq:ngest10}
\begin{split}
\mathcal{E}_{g}(f_p(\l)(\s))
\leq &
\max_{i=1, \cdots, p}\mathcal{E}_{g}(\varphi_{i})\frac
{
\sum_{i=1}^p\alpha_{i}^{2}\mathcal{D}_{g}(\varphi_{i})
}
{
\left(\sum_{i=1}^p\left(\alpha_{i}^{2}\mathcal{D}_{g}(\varphi_{i})\right)^{\frac{n}{n-2}}\right)^{\frac{n-2}{n}}
}
\\ 
& +
\frac
{\left(1+o_{\sum_{i\neq j}\varepsilon_{i, j}}(1)\right)\mathcal {Y}(S^n_+)}
{\left(\sum_{i=1}^p\alpha_{i}^{\frac{2n}{n-2}}\right)^{\frac{n-2}{n}}}\left(
\sum_{i\neq j}\left[1-2\frac{\alpha_{i}^{\frac{4}{n-2}}\left(\sum_{i=1}^p\alpha_{i}^{2}\right)}{\left(\sum_{i=1}^p\alpha_{i}^{\frac{2n}{n-2}}\right)}\right]\alpha_{i}\alpha_{j}\frac{\varepsilon_{i,j}}{c_1}\right).
\end{split}\end{equation} 
This inequality has the following impact. First note, that the function
\begin{equation}\label{eq:defgamma}
\begin{split} 
\Gamma:\{\gamma \in [0,1]^{p}:&\;\;\;  \sum_{i=1}^p\gamma_{i} =1\}\longrightarrow  \mathbb{R}_{+}\\&\gamma
\longrightarrow 
\frac{\sum_{i=1}^p\gamma_{i}^{2}}{(\sum_{i=1}^p\gamma_{i}^{\frac{2n}{n-2}})^{\frac{n-2}{n}}}
\end{split}\end{equation} 
has the strict global maximum 
\begin{equation}\label{eq:gammamax}
\gamma_{\max}=\left(\frac{1}{p},\ldots,\frac{1}{p}\right)
\end{equation} 
with \;$\Gamma(\gamma_{\max})=p^{\frac{2}{n}}$. Thus, using Lemma \ref{eq:bsvolestl}, Lemma \ref{eq:bubbleestl}, \eqref{eq:numden}, \eqref{eq:varepijep}, and \eqref{eq:ngest10},
we infer that for any \;$\nu>0$, for every \;$\epsilon>0$\;and small, and for every \;$p\in\N^*$, there exists \;$\l_p:=\l_p(\nu, \epsilon)\geq\frac{2}{\d_0^{\frac{n}{2}}}$\; such for every \;$\l\geq \l_p$\; and for every $\sigma:=\sum_{i=1}^p\alpha_i\d_{a_i}\in B_p(\partial M)$, there holds
\begin{equation}\begin{split} 
\mathcal{E}_{g}(f_p(\l)(\s))<p^{\frac{2}{n}}\mathcal{Y}(S^n_+),
\end{split}\end{equation} 
whenever 
\begin{equation}\begin{split} 
\exists i_o\neq j_0 \;\;\;\text{such that}\;\;\; \frac{\alpha_{i_0}}{\alpha_{j_0}}>\nu\;
\;\;\;\; \text{and}\;\;\;\; 
\sum_{i\neq j}\varepsilon_{i,j}\leq\epsilon,
\end{split}\end{equation} 
thereby ending the proof of point 2). Now, we are going to treat point 3) and end the proof of the Lemma. Thus, we may assume
\begin{equation}\label{eq:assump1}
\begin{split} 
\forall\;\, i,j\;\;\frac{\alpha_{i}}{\alpha_{j}}= 1+o^+_{\mu}(1)\;\;\;\;
\;\text{and}\;\;\; \sum_{i\neq j}\epsilon_{i,j}\ll 1,
\end{split}\end{equation} 
where \;$o_{\mu}^+(1)$\; is a positive quantity depending only \;$\mu$\; with \;$\mu$\; small and verifying the property that it tends to $0$ as \;$\mu$\; tends to \;$0$. So, using \eqref{eq:ngest10}, \eqref{eq:assump1}, and the properties of $\Gamma$ (see \eqref{eq:defgamma} and \eqref{eq:gammamax}), we infer that the following estimate holds
\begin{equation}\label{eq:ngest11}
\begin{split} 
\mathcal{E}_{g}(f_p(\l)(\s))
\leq &
\max_{i=1, \cdots, p}\mathcal{E}_{g}(\varphi_{i})p^{\frac{2}{n}}\\&
-(1+o_{\sum_{i\neq j}\varepsilon_{i, j}}(1)+o_{\mu}(1))
\frac
{\mathcal{Y}(S^n_+)}
{(\sum_{i=1}^p\alpha_{i}^{\frac{2n}{n-2}})^{\frac{n-2}{n}}}
\sum_{i\neq j}\alpha_{i}\alpha_{j}\frac{\varepsilon_{i,j}}{c_1} \\
= &
\max_{i=1, \cdots, p}\mathcal{E}_{g}(\varphi_{i})p^{\frac{2}{n}}
-(1+o_{\sum_{i\neq j}\varepsilon_{i, j}}(1)+o_{\mu}(1))
\mathcal{Y}(S^n_+)p^{\frac{2-n}{n}}
\sum_{i\neq j}\frac{\varepsilon_{i,j}}{c_1}.
\end{split}\end{equation} 
Now, using Lemma \ref{eq:bubbleestl}, \eqref{eq:varepijep},  \eqref{eq:assump1}, and \eqref{eq:ngest11}, we have that there exists \;$C_0>0$, $\nu_0>1$, $\l_0\geq \frac{2}{\delta_0^{\frac{n}{2}}}$\; and \;$0<\epsilon_0\leq\delta_0$\; such that for every \;$1<\nu\leq \nu_0$, for every \;$0<\epsilon\leq\epsilon_0$, for every \;$p\in \N^*$, for every \;$\l\geq \l_0$, and for every \;$\sigma:=\sum_{i=1}^p\alpha_i\d_{a_i}\in B_p(\partial M)$, we have if $\frac{\alpha_{i}}{\alpha_{j}} \leq\nu$\;\;\;$\forall i, j$\;\;
\; and\;\;\; $\sum_{i\neq j}\epsilon_{i,j}\leq\epsilon$, then there holds
\begin{equation}\label{eq:nges12}
\begin{split} 
\mathcal{E}_{g}(f_p(\l)(\s))
\leq &
p^{\frac{2}{n}}\mathcal{Y}(S^n_+)
\left(
1+\frac{C_0}{\lambda^{n-2}}
-
\frac{1}{2c_1p}
\sum_{i\neq j}\varepsilon_{i,j}
\right).
\end{split}\end{equation} 
Thus, recalling that (see \eqref{eq:48}) 
\begin{equation}\begin{split}
\varepsilon_{i,j} = (1+o_{\varepsilon_{i, j}}(1))c_{3}\frac{G(a_{i},a_{j})}{\lambda^{n-2}},
\end{split}\end{equation}
and using again \eqref{eq:varepijep}, we infer from \eqref{eq:nges12} that up to taking $\epsilon_0$ smaller, for every $1<\nu\leq \nu_0$, for every $0<\epsilon\leq\epsilon_0$, for every \;$p\in \N^*$, for every $\l\geq \l_0$, and for every $\sigma:=\sum_{i=1}^p\alpha_i\d_{a_i}\in B_p(\partial M)$, there holds  $$\frac{\alpha_{i}}{\alpha_{j}} \leq\nu\;\;\;\forall i, j\;\;
\; \text{and}\;\;\; \sum_{i\neq j}\epsilon_{i,j}\leq\epsilon$$ imply 
\begin{equation}\begin{split} 
\mathcal{E}_{g}(f_p(\l)(\s))
\leq &
p^{\frac{2}{n}}\mathcal{Y}(S^n_+)\left
(
1+\frac{C_0}{\lambda^{n-2}}
-
\frac{c_{3}}{4c_1p\lambda^{n-2}}
\sum_{i\neq j}G(a_{i},a_{j})
\right)\\
\leq &
p^{\frac{2}{n}}\mathcal{Y}(S^n_+)\left(
1+\frac{C_0}{\lambda^{n-2}}
-
c_g\frac{(p-1)}{\lambda^{n-2}}
\right),
\end{split}\end{equation} 
thereby ending the proof of point 3), and hence of the Lemma.
\end{pf}
\vspace{4pt}

\noindent
\begin{pfn}{ of Proposition \ref{eq:baryest}}\\
It follows from Lemma \ref{eq:baryestaux} by taking \;$C_0$\; and \;$\nu_0$\; to be the ones given by Lemma \ref{eq:baryestaux}, $\varepsilon_0:=\frac{\epsilon_0}{2}$, and \;$\l_p:=\l_p(\varepsilon, \nu_0):=\max\{\l_p(\frac{\varepsilon}{2}), \l_p( 2\varepsilon, \nu_0), \l_0\}$, where\; $\epsilon_0$, $\l_p(\frac{\varepsilon}{2})$, $\l_p(2\varepsilon, \nu_0)$, and \;$\l_0$\; are as in Lemma \ref{eq:baryestaux}.
\end{pfn}
\vspace{6pt}

\noindent
Now, using Proposition \ref{eq:baryest}, we have the following corollary which provides the second and decisive step in our goal of showing that the Brendle-Chen\cite{bs}'s bubbles are algebraic topological ones. Indeed, the latter corollary will be used together with Proposition \ref{eq:baryest} in the next section  to carry a suitable scheme of the algebraic topological argument of Bahri-Coron\cite{bc}.
\begin{cor}\label{eq:largep}
There exists $p_0\in \N^*$ large enough such that for every \;$0<\varepsilon\leq\varepsilon_0$, and for every for every $\l\geq\l_{p_0}$ (where \;$\varepsilon_0$\; and \;$\l_{p_0}$\; are given by Proposition \ref{eq:baryest}), there holds 
$$
\mathcal{E}_g(f_{p_0}(\l))(B_{p_0}(\partial M)))\subset W_{p_0-1}.
$$
\end{cor}
\begin{pf}
It follows directly from Proposition \ref{eq:baryest} and the definition of \;$W_{p_0-1}$ (see \eqref{eq:defenergylevel} with \;$p$\; replaced by \;$p_0-1$).
\end{pf}

\section{Algebraic topological argument for existence}\label{eq:algtop}
In this section, we show that the Brendle-Chen\cite{bs}'s bubbles are algebraic topological ones by using directly the latter section to carry a suitable scheme of the barycenter technique of Bahri-Coron\cite{bc}. To do so, we first introduce the notion of {\em neighborhood of potential critical points at infinity} of the Euler-Lagrange functional \;$\mathcal{E}_g$. Precisely, for $p\in \N^*$, $0<\varepsilon\leq \varepsilon_0$ (where $\varepsilon_0$ is given by Proposition \ref{eq:baryest}), and \;$(l, q)\in \N^2$ with $2l+q=p$, we define \;$V(l, q, \varepsilon)$\; the \;$(l, q, \varepsilon)$-neighborhood of potential critical points at infinity of $\mathcal{E}_g$, namely
\begin{equation}\label{eq:vqs}
\begin{split}
V(l, q, \varepsilon):=\{u\in W^{1, 2}_+(\ov M):\;\;\exists a_1, \cdots, a_l\in M, \;\;a_{l+1}, \cdots, a_{l+q}\in \partial M,\;\;\alpha_1, \cdots, \alpha_{l+q}>0, \;\l_1, \cdots,\l_{l+q}\geq\frac{1}{\varepsilon},\\||u-\sum_{i=1}^{l+q}\alpha_i\varphi_{a_i, \l_i}||\leq \varepsilon, \;\;\frac{\alpha_i}{\alpha_j}\leq \nu_0\;\;\;\;\text{and}\;\;\;\varepsilon_{i, j}\leq \varepsilon, \;\;i\neq j=1, \cdots, l+q\},
\end{split}
\end{equation}
where $||\cdot||$\; denotes the standard \;$W^{1, 2}$-norm, $\varepsilon_{i, j}:=\varepsilon_{i, j}(A, \bar \l)$ with $A:=(a_1, \cdots, a_{l+q})$, $\bar \l:=(\l_1, \cdots, \l_{l+q})$, \; $(\varepsilon_{i, j}(A, \bar \l))$'s are defined by \eqref{eq:varepsilonij}, $\varphi_{a_i, \l}$\; is given by \eqref{eq:varphialb} for $i=l+1, \cdots, l+q$, $\varphi_{a_i, \l}$ \;is given by \eqref{eq:varphiali} for \;$i=1, \cdots, l$, and\; $\nu_0$ is given by Lemma \ref{eq:baryest}. 
Using the $(l, q, \varepsilon)$-neighborhoods of potential critical points at infinity of \;$\mathcal{E}_g$, we define\; $V(p, \varepsilon)$\;the \;$(p, \varepsilon)$-neighborhood of potential critical points at infinity of \;$\mathcal{E}_g$\; as follows
\begin{equation}
V(p, \varepsilon)=\cup_{(l, q):\; 2l+q=p}\;V(l, q, \varepsilon).
\end{equation}
\vspace{6pt}

\noindent
Concerning the sets $V(0, p, \varepsilon)$, for every $p\in \N^*$, we have that there exists $0<\varepsilon_p\leq\varepsilon_0$ (where $\varepsilon_0$ is given by Lemma \ref{eq:baryest}) such that for every \;$0<\varepsilon\leq \varepsilon_p$, we have that
\begin{equation}\label{eq:mini}
\begin{split}
\forall u\in V(0, p, \varepsilon),\;\;\; \text{the minimization problem }\;\;\;\min_{B_{C_1\varepsilon}^{p}}||u-\sum_{i=1}^{p}\alpha_i\varphi_{a_i, \l_i}||
\end{split}
\end{equation}
has a unique solution, up to permutations, where $B^{p}_{C_1\varepsilon}$ is defined as follows
\begin{equation}
\begin{split}
B_{C_1\varepsilon}^{p}:=\{(\bar\alpha, A, \bar \l)\in \R^{p}_+\times (\partial M)^p\times (0, +\infty)^{p}:\;\;\l_i\geq \frac{1}{\epsilon}, i=1, \cdots, p,\\\frac{\alpha_i}{\alpha_j}\leq \nu_0\;\; \text{and}\;\ \varepsilon_{i, j}\leq C_1\varepsilon, i\neq j=1, \cdots, p\},
\end{split}
\end{equation}
with $C_1>1$. Furthermore, we define the selection map $s_{p}: V(0, p, \varepsilon)\longrightarrow (\partial M)^p/\sigma_p$ as follows
\begin{equation}\label{eq:select}
s_{p}(u):=A, \;\;u\in V(0, p, \varepsilon), \;\,\text{and} \,\;A\;\;\text{is given by}\;\,\eqref{eq:mini},
\end{equation}
\begin{rem}\label{eq:remselect}
We would like to emphasize the fact that we have defined a selection map only on \;$V(0, p, \varepsilon)$ and not on the whole \;$V(p, \varepsilon)$. This is due to the fact that, our algebraic topological argument uses only bubbles centered at boundary points and that a continuous extension on \;$V(p, \varepsilon)$\; can also be achieved if the values it takes on \;$V(p, \varepsilon)\setminus V(0, p, \varepsilon)$ do not matter, which is clearly the case here. However, using the Djadli-Malchiodi\cite{dm} barycentric projection map, precisely its boundary adaptation of Ndiaye\cite{nd2}, one can even achieve a continuous extension with the property of {\em homotopic to identity} in the sense of Djadli-Malchiodi\cite{dm}.
\end{rem}

\vspace{6pt}

\noindent
Now having introduced the neighborhoods of potential critical points at infinity of the Euler-Lagrange functional \;$\mathcal{E}_g$, we are ready to present our algebraic topological argument for existence. In order to do that, we start by the following classical deformation Lemma which follows from the same arguments as for its counterparts in classical application of the algebraic topological argument of Bahri-Coron\cite{bc}(see for example Proposition 6 in \cite{bc} or Lemma 17 in \cite{BB}) and the fact that the \;$\varphi_{a, \l}$ can replace the standard bubbles in the analysis of diverging PS sequences of the Euler-Lagrange functional\; $\mathcal{E}_g$.
\begin{lem}\label{eq:classicdeform}
Assuming that \;$\mathcal{E}_g$ \;has no critical points, then for every \;$p\in \N^*$, up to taking $\varepsilon_p$\; smaller (where \;$\varepsilon_p$\; is given by \eqref{eq:mini}), we have that for every\; $0<\varepsilon\leq \varepsilon_p$, there holds \;$(W_p,\; W_{p-1})$\; retracts by deformation onto \;$(W_{p-1}\cup A_p, \;W_{p-1})$\; with \;$V(p, \;\tilde \varepsilon)\subset A_p\subset V(p, \;\varepsilon)$\; where \;$0<\tilde \varepsilon<\frac{\varepsilon}{4}$\; is a very small positive real number and depends on $\varepsilon$.
\end{lem}
\vspace{6pt}

\noindent
Using Proposition \ref{eq:baryest} and Lemma \ref{eq:classicdeform}, we are going to show that if \; $\mathcal{E}_g$\; has no critical points, then for \;$\l$\; large enough, the map \;$(f_1(\l))_*$\; is  well defined and not trivial in \;$H_{n-1}(W_1,\; W_0)$. Precisely, we show:
\begin{lem}\label{eq:nontrivialf1}
Assuming that \;$\mathcal{E}_g$\; has no critical points and \:$0<\varepsilon\leq  \varepsilon_1$ (where \;$\varepsilon_1$\; is given by \eqref{eq:mini}), then up to taking \;$\varepsilon_1$\; smaller and \;$\l_1$\; larger (where \;$\l_1$\; is given by Proposition \ref{eq:baryest}), we have that for every \;$\l\geq \l_1$, there holds
$$
f_1(\l): \;(B_1(\partial M),\; B_0(\partial M))\longrightarrow (W_1, \;W_0)
$$
is well defined and satisfies
$$
(f_1(\l))_*(w_1)\neq 0\;\;\;\;\text{in}\;\;\;\;H_{n-1}(W_1, \;W_0).
$$
\end{lem}
\begin{pf}
It follows from the selection map \;$s_{1}$\; given by \eqref{eq:select}, Proposition \ref{eq:baryest} and the same arguments as in Lemma 26 in \cite{gam2}.
\end{pf}
\vspace{6pt}

\noindent
Next, using Lemma \ref{eq:transfert}, Proposition \ref{eq:baryest}, Lemma \ref{eq:classicdeform}, and the algebraic topological argument of Bahri-Coron\cite{bc} in the form developed by the second author\cite{nd6} for the study of the supercritical \;$Q$-curvature problem, we are going to show that if for \;$\l$ large the orientation class \;$w_p$\; survives ``topologically`` the embedding into \;$(W_p, \;W_{p-1})$ via $f_p(\l)$, then for $\l$ large \;$w_{p+1}$\; survives ''topologically`` the embedding into \;$(W_{p+1},\; W_p)$ via $f_{p+1}(\l)$. Precisely, we prove the following proposition:
\begin{pro}\label{eq:nontrivialrecursive}
Assuming that \;$\mathcal{E}_g$\; has no critical points and \;$0<\varepsilon\leq  \varepsilon_{p+1}$ (where \;$\varepsilon_{p+1}$\; is given by \eqref{eq:mini}), then up to taking \;$\varepsilon_{p+1}$\; smaller, and \;$\l_p$\; and \;$\l_{p+1}$\; larger\;(where \;$\l_p$\; and \;$\l_{p+1}$\; are given by Proposition \ref{eq:baryest}), we have that for every \;$\l\geq \max\{\l_p, \l_{p+1}\}$, there holds
$$
f_{p+1}(\l): (B_{p+1}(\partial M),\; B_{p}(\partial M))\longrightarrow (W_{p+1}, \;W_{p})
$$
and 
$$
f_p(\l): (B_p(\partial M), \;B_{p-1}(\partial M))\longrightarrow (W_p, \; W_{p-1})
$$
are well defined and satisfy
$$(f_p(\l))_*(w_p)\neq 0\;\;\;\; \text{in}\;\; \;\;H_{np-1}(W_p, \;W_{p-1})$$ implies
$$(f_{p+1}(\l))_*(w_{p+1})\neq 0\;\;\;\; \text{in} \;\;\;\;H_{n(p+1)-1}(W_{p+1}, \;W_{p}).$$
\end{pro}
\begin{pf}
First of all, we let \;$p\in \N^*$\; and \;$0<\varepsilon_{p+1}$, where \;$\varepsilon_{p+1}$\; is given by Proposition \ref{eq:baryest}. Next, recalling that we have assumed that \;$\mathcal{E}_g$\; has no critical points, and using Lemma \ref{eq:classicdeform}, then up to taking \;$\varepsilon_{p+1}$\; smaller, we infer that the following holds
\begin{equation}\label{eq:idents4}
(W_{p+1}, \;W_{p})\simeq (W_{p}\cup \mathcal{A}_{p+1}, \;W_{p}),
\end{equation}
with
\begin{equation}\label{eq:infinitynontrivial}
V(0, \;p+1, \;\tilde\varepsilon)\subset V(p+1, \;\tilde \varepsilon)\subset \mathcal{A}_{p+1}\subset V(p+1,\; \varepsilon), \;\;0<4\tilde\varepsilon<\varepsilon.
\end{equation}
Now, using Lemma \ref{eq:interactestl} and Proposition \ref{eq:baryest}, we have that for every \;$\l\geq \max\{\l_p, \l_{p+1}\}$ (where $\l_p$ and $\l_{p+1}$\; are given by Proposition \ref{eq:baryest}), there holds
\begin{equation}\label{eq:mappingfp1}
f_{p+1}(\l): (B_{p+1}(\partial M), \;B_{p}(\partial M))\longrightarrow (W_{p+1}, \;W_{p}),
\end{equation}
and
\begin{equation}\label{eq:mappingfp}
f_p(\l): (B_p(\partial M), \; B_{p-1}(\partial M))\longrightarrow (W_p, \;W_{p-1}),
\end{equation}
are well defined and hence have that the first point is proven. Next, using Proposition \ref{eq:baryest}, \eqref{eq:mappingfp1}, and \eqref{eq:mappingfp}, we have that up to taking \;$\l_{p+1}$\; and \;$\l_{p}$\; larger (for example larger than \;$4\max\{\l_{p+1}(\varepsilon),\;\l_{p}(\varepsilon), \;\l_{p}(2\tilde\varepsilon), \l_p(\frac{\tilde\varepsilon}{2}), \frac{1}{\tilde \varepsilon}\}$, where $\l_p(\varepsilon)$, \;$\l_{p+1}(\varepsilon), \l_p(2\tilde\varepsilon)$,\;and \;$\l_{p}(\frac{\tilde\varepsilon}{2})$\; are given by Proposition \ref{eq:baryest} and \;$\tilde \varepsilon$ \; is given by \eqref {eq:infinitynontrivial}) the following diagram 
\begin{equation}\label{eq:diag2}
\begin{CD}
 (B_{p+1 }(\partial M), \;\mathcal{O}(B_{p}(\partial M)))  @>f_{p+1}(\lambda)>> (W_{p+1}, \; W_{p})\\
@AAA @AAA\\
(\mathcal{O}(B_{p}(\partial M)),\;B_{p-1}(\partial M))  @>f_{p}(\l)>> (W_{p}, \; W_{p-1})
\end{CD}
\end{equation}
is well defined and commutes, where
\begin{equation}\label{eq:nlbary}
\mathcal{O}(B_{p}(\partial M)):=\{\sigma=\sum_{i=1}^{p+1}\alpha_i\d_{a_i}\in B_{p+1}(\partial M):\;\;\exists \;i_0\neq j_0:\;\frac{\alpha_{i_0}}{\alpha_{j_0}}>\nu_0\;\;\;\;\text{or}\;\;\;\;\sum_{i\neq j}\varepsilon_{i, j}> \tilde\varepsilon\},
\end{equation}
with \;$\nu_0$\; given by Proposition \ref{eq:baryest}. On the other hand, we have
\begin{equation}\label{eq:isodel1}
B_{p+1}(\partial M)\setminus \mathcal{O}(B_{p}(\partial M))\simeq B_{p+1}(\partial M)\setminus B_{p}(\partial M),
\end{equation}
and
\begin{equation}\label{eq:isodel2}
\mathcal{O}(B_{p}(\partial M))\simeq  B_{p}(\partial M).
\end{equation}
Now, using \eqref{eq:purem}, Lemma \ref{eq:transfert}, and \eqref{eq:isodel1}, we derive
\begin{equation}\label{eq:tranfert2}
\begin{CD}
 H^{n-1}(B_{p+1}(\partial M)\setminus \mathcal{O}(B_{p}(\partial M)))\times H_{n(p+1)-1}(B_{p+1}(\partial M), \;B_{p}(\partial M))@>\frown>> H_{n(p+1)-n}(B_{p+1}(\partial M),\; B_{p}(\partial M))\\@>\partial>>H_{n(p+1)-n-1}(B_{p}(\partial M),\; B_{p-1}(\partial M)).
 \end{CD}
\end{equation}
Furthermore, using \eqref{eq:idents4}, we infer that
\begin{equation}\label{eq:transfert3}
\begin{CD}
 H^{n-1}(\mathcal{A}_{p+1})\times H_{n(p+1)-1}(W_{p+1},\; W_{p})@>\frown>> H_{n(p+1)-n}(
 W_{p+1},\; W_{p})\\@>\partial>>H_{n(p+1)-n-1}(W_{p},\; W_{p}).
 \end{CD}
\end{equation}
Moreover, passing to homologies in \eqref{eq:diag2} and using \eqref{eq:isodel2}, we derive that
\begin{equation}\label{eq:diag3a}
\begin{CD}
 H_{n(p+1)-1}(B_{p+1}(\partial M), \;B_{p}(\partial M))   &@>(f_{p+1}(\l))_*>>  H_{n(p+1)-1}(W_{p+1}, \; W_{p}),
\end{CD}
\end{equation}
and
\begin{equation}\label{eq:diag3b}
\begin{CD}
 H_{np-1}(B_{p}(\partial M), \;B_{p-1}(\partial M))  & @>(f_{p}(\l))_*>>  H_{np-1}(W_{p}, \;W_{p-1})
\end{CD}
\end{equation}
are well defined and the following diagram commutes
\begin{equation}\label{eq:diag3}
\begin{CD}
 H_{n(p+1)-n}(B_{p+1}(\partial M), \;B_{p}(\partial M))   &@>(f_{p+1}(\l))_*>>  H_{n(p+1)-n}(W_{p+1},\; W_{p})\\
@V {\partial} VV& @ V {\partial} VV\\
H_{np-1}(B_{p}(\partial M), \;B_{p-1}(\partial M))  & @>(f_{p}(\l))_*>>  H_{np-1}(W_{p}, \;W_{p-1}).
\end{CD}
\end{equation}
Next, recalling that we have taken \;$\l_{p+1}$\; and \;$\l_{p}$\; larger  than \;$4\max\{\l_{p+1}(\varepsilon),\;\l_{p}(\varepsilon), \;\l_{p}(2\tilde \varepsilon), \l_p(\frac{\tilde\varepsilon}{2}), \frac{1}{\tilde \varepsilon}\}$, we derive that
\begin{equation}\label{eq:insidevm}
f_{p+1}(\l)\left(B_{p+1}(\partial M)\setminus \mathcal{O}(B_{p}(\partial M))\right)\subset V(0, \;p+1,\;\tilde\varepsilon)\subset V(p+1, \;\tilde \varepsilon)\subset \mathcal{A}_{p+1}\subset V(p+1, \varepsilon).
\end{equation}
Thus, using \eqref{eq:defom1}, \eqref{eq:select}, \eqref{eq:insidevm}, and recalling Remark \ref{eq:remselect}, we infer that 
\begin{equation}\label{eq:rbc}
(f_{p+1}(\l))^*(s_{p+1}^*(O^*_{\partial M}))=O^*_{\partial M}\;\;\;\text{with}\;\; \;s^*_{p+1}(O^*_{\partial M})\neq 0\;\;\;\text{in}\;\;\; H^{n-1}(\mathcal{A}_{p+1}).
\end{equation}
On the other hand, using \eqref{eq:diag3}, we derive that
\begin{equation}\label{eq:commutation}
\partial (f_{p+1}(\l))_{*}  =  (f_{p}(\l))_{*}\partial \;\;\;\text{in}\;\;\;H_{n(p+1)-n}(B_{p+1}(M), \;B_{p}(\partial M)).
\end{equation}
Now, combining \eqref{eq:purem}, Lemma \ref{eq:transfert}, \eqref{eq:isodel1}, \eqref{eq:tranfert2}, \eqref{eq:transfert3}, \eqref{eq:diag3a}, \eqref{eq:rbc}, and \eqref{eq:commutation}, we obtain
\begin{equation}\label{eq:transfim}
\begin{split}
(f_{p}(\l))_*(\o_{p})&=(f_{p}(\l))_*\left(\partial(O^{*}_{\partial M} \smallfrown \o_{p+1}) \right)\\
&=(f_{p}(\l))_*\left(\partial(((f_{p+1}(\l))^*(s_{p+1}^{*}(O^{*}_{\partial M}))) \smallfrown \o_{p+1}) \right)\\
&=\partial\left((f_{p+1}(\l))_*(((f_{p+1}(\l))^*(s_{p+1}^{*}(O^{*}_{\partial M})) \smallfrown \o_{p+1}) \right)\\
&=\partial(s_{p+1}^{*}(O^{*}_{\partial M}) \smallfrown ((f_{p+1}(\l))_*(\o_{p+1}))),
\end{split}
\end{equation}
with all the equalities holding in \;$H_{np-1}(W_{p}, \;W_{p-1})$. Hence, clearly, \eqref{eq:transfim} and the assumption 
$$(f_p(\l))_*(w_p)\neq 0\;\; \;\text{in}\; \;\;H_{np-1}(W_p, \;W_{p-1})$$ implies
$$(f_{p+1}(\l))_*(w_{p+1})\neq 0\;\;\; \text{in} \;\;\;H_{n(p+1)-1}(W_{p+1}, \;W_{p}),$$
as desired, thereby completing the proof of Proposition \ref{eq:nontrivialrecursive}.
\end{pf}
\vspace{4pt}

\noindent
Now, we are ready to present the proof of Theorem \ref{eq:existence}.\\\\
\begin{pfn}{ of Theorem \ref{eq:existence}}\\
It follows by a contradiction argument from Corollary \ref{eq:largep}, Lemma \ref{eq:nontrivialf1} and Proposition \ref{eq:nontrivialrecursive}.
\end{pfn}

\end{document}